\documentclass[12pt,reqno]{amsart}
\usepackage{amssymb,amsmath,amsthm,longtable,multirow,array}

\usepackage{color}
\usepackage{mathrsfs}
\numberwithin{equation}{section}

\newtheorem{Theorem}{Theorem}[section]
\newtheorem{Proposition}[Theorem]{Proposition}

\newtheorem{Lemma}[Theorem]{Lemma}
\newtheorem{Corollary}[Theorem]{Corollary}
\theoremstyle{definition}
\newtheorem{Definition}[Theorem]{Definition}
\theoremstyle{remark}
\newtheorem{Example}[Theorem]{Example}

\newtheorem*{remark*}{Remark}

 \providecommand\CC{{\mathbb C}}

\providecommand\RR{{\mathbb R}}

 \providecommand\ZZ{{\mathbb Z}}
\providecommand\HH{\mathbb H}

\providecommand\brs{\begin{remark*}}
\providecommand\ers{\end{remark*}}

\providecommand\om{\omega}

\providecommand\lm{\lambda}

\providecommand\be{\begin{enumerate}}
\providecommand\ee{\end{enumerate}}
\providecommand\bT{\begin{Theorem}}
\providecommand\eT{\end{Theorem}}
\providecommand\bP{\begin{Proposition}}
\providecommand\eP{\end{Proposition}}
\providecommand\bD{\begin{Definition}}
\providecommand\eD{\end{Definition}}
\providecommand\bE{\begin{Example}}
\providecommand\eE{\end{Example}}
 \providecommand\bL{\begin{Lemma}}
\providecommand\eL{\end{Lemma}}
\providecommand\bC{\begin{Corollary}}
\providecommand\eC{\end{Corollary}}

\providecommand\bpp{\begin{proof}} \providecommand\epp{\end{proof}}
\providecommand\bee{\begin{equation}}
\providecommand\eee{\end{equation}}

  \providecommand\beqq{\begin{eqnarray*}}
\providecommand\eeqq{\end{eqnarray*}}
\providecommand\a{\alpha}

\providecommand\th{\theta}
 
\providecommand\ol{\overline}
 \providecommand\rt{\rightarrow}
\providecommand\?{\infty}
 \providecommand\ul{\underline}

\providecommand\besp{\begin{split}}
\providecommand\eesp{\end{split}}
\providecommand\bay{\begin{array}}
\providecommand\eay{\end{array}}
\providecommand\sgn{\operatorname{sgn}}
\providecommand\bes{\begin{equation}\begin{split}}
\def\a{\alpha}
\def\b{\beta}

\def\f#1#2{\frac{#1}{#2}}
\def\t{\tau}
\def\th{\theta}
\def\s{\sigma}
\def\z{\zeta}

\def\<{\leq}
\begin{document}

\title[Hilbert Transformation]{Hilbert Transformation and  $r\mathrm{Spin}(n)+\RR^n$  Group}
\author{Pei Dang}
\author{Hua Liu}
\author{Tao Qian}

\begin{abstract} In this paper we study symmetry properties of the Hilbert transformation of several real variables in the Clifford algebra setting. In order to describe the symmetry properties we introduce the group $r\mathrm{Spin}(n)+\RR^n, r>0,$ which is essentially an extension of the ax+b group.
The study concludes that the Hilbert transformation has certain characteristic symmetry properties in terms of $r\mathrm{Spin}(n)+\RR^n.$  In the present paper, for $n=2$ and $3$ we obtain, explicitly,
the induced spinor representations of the
$r\mathrm{Spin}(n)+\RR^n$ group.
Then we decompose the natural representation of $r\mathrm{Spin}(n)+\RR^n$  into the direct sum of some two irreducible spinor representations, by which we characterize the Hilbert transformation in $\RR^3$ and $\RR^2.$ Precisely, we show that a nontrivial skew operator is the Hilbert transformation if and only if it is invariant under the action of the  $r\mathrm{Spin}(n)+\RR^n, n=2,3,$ group.

\end{abstract}

\address{pdang@must.edu.mo,  Faculty of Information Technology, Macau University of Science and Technology,
Macau }
\address{hualiu@tute.edu.cn, Department of Mathematics, Tianjin
University of Technology and Education, Tianjin 300222, China}
\address{fsttq@umac.mo, Department of Mathematics, Macau University, Macau }

\maketitle

\bigskip

\section{Introduction}

The underlying function space of this study is $L_2(\RR^n).$ The one dimensional Hilbert transformation is defined as
\begin{equation}\label{hilbert1}
\nonumber (H^0f)(x)\triangleq\f1\pi\ \lim_{\epsilon\to 0}\int_{|x-y|>\epsilon}\frac{f(y)}{x-y}dy=\f1\pi\ {\rm p.v.} \int_{-\infty}^\infty \frac{f(y)}{x-y}dy,\ x\in\RR,
\end{equation}
where the integral is considered as the extension to $L_2(\RR^n)$ of the Cauchy
principle-valued integral operator over the subspace of $L_2(\RR)$ consisting of the functions of the H\"{o}lder continuity.  The extension is based on the $L_2$-boundedness of the Hilbert transformation operator, \cite{stein}. Hilbert transformations in various contexts have an ample amount of applications, including solving problems in aerodynamics, condensed matter physics, optics, fluids, and engineering (see, for instance, \cite{king}). In particular, the Hilbert transformation plays an
important role in analytic boundary value problems, harmonic and signal analysis.

 There are different profiles of the operator in  diverted  fields of mathematics, of which one is its symmetry. In \cite{king} it is pointed out that the one dimensional Hilbert transformation is invariant under the actions of  the translation and the dilation. The translation and the dilation generate the \bf ax+b \rm group,  which is the
group of the affine transformations $x\rt ax+b$ of $\RR$ with $a>0$
and $b\in \RR$. Its underlying manifold is $(0,\?)\times\RR.$ The
group law is defined by
$$
((a,b)(a',b'))(x)=aa'x+b+ab'=(aa',b+ab')(x),\ x\in\RR,
$$
which gives
\begin{equation}\label{mult}
\nonumber   (a,b)(a',b')=(aa',b+ab').
\end{equation}
It is  easy to check that $(a,b)^{-1}=(\f1a,-\f ba)$, $da/a^2$ is the
left-Haar measure and $dadb/a$ is the right-Haar measure of this
group\ \cite{folland}.

There exists a natural unitary infinite
dimension representation $\pi$ of $G$ over the Hilbert space $L_2(\RR)$. Denote by
$\mathfrak{U}(L_2(\RR))$ the operator group of the unitary
auto-isomorphism of $L_2(\RR)$. Then the group morphism
$\pi:G\rightarrow\mathfrak{U}(L_2(\RR))$ is defined by
\begin{equation}\label{nature1}
 \nonumber   (\pi(a,b)f)(x)=(\f1a)^{\f12}f(\f{x-b}a),\ x\in\RR.
\end{equation}
$\pi(a,b)$ is also written as $\pi_{ab}$.

The Hilbert transformation is invariant under the action of the representation, i.e.,  $\pi$ commutes with $H$,
\begin{equation}\label{com1}
    (H^0\pi_{ab})f=(\pi_{ab}H^0)f,\ \quad (a,b)\in G,\ f\in L_2(\RR).
\end{equation}
In \cite{eh}, A.S. Elmabrok and  O. Hutnik considered wavelet transformations  from the viewpoint of  induced representation of this group. It shows that the ax+b group has close relevance with physics.  In fact, using Gelfand's representation of ax+b, the Hilbert
transformation may be, in depth, characterized by the commutativity (\ref{com1}), \cite{dlq}.
\vskip 3mm

 In \cite{st}, E.M. Stein refers to the following variation of the Hilbert transformation:
\begin{equation}\label{hilbert-vec}
  H_vf(x)=\f1\pi\int_{-\epsilon}^\epsilon f(x-v(x)y)\f{dy}y,
\end{equation}
where $v(x)$ is a Lipschitz function and $H_v$ is called the directional Hilbert transformation. There are other formulations as well (\cite{kim}, \cite{ch}, \cite{CMcM}).
Among several equivalent definitions of the Hilbert transformation the one directly related to
 boundary limits of analytic functions, or precisely the Plemelj formula, would be the most
 original and with the widest connections to applications. In general, with the one complex variable
 setting, the Hilbert transformation on a manifold is defined to be the mapping from the real part to the
 imaginary part of the boundary limit of an analytic function on one of the two regions cut out by the manifold ([4]).
 Such definition can be generalized to higher dimensions with the Cauchy structure of several complex variables
 or alternatively  that of the Clifford algebra.
  In  $\RR^n$, some researchers define the Hilbert transformation to be (see, for instance, \cite{king})
\begin{equation}\label{hilbert-king}
  H_nf(x)=\f1 {\pi^n}\prod_{j=1}^n\lim_{\varepsilon_j\rt0}\int_{|x_j-s_j|>\varepsilon_j}f(s)\prod_{k=1}^n\f{ds_k}{x_k-s_k}.
\end{equation}
$H_n$ can be written as a tensor product of the partial Hilbert transformations, i.e., $H_n=\prod_{k=1}^nH_{(k)}$, where
\begin{equation}\label{hilbert-king1}
  H_{(k)}f(s_1,\cdots,s_{k-1},x_k,s_{(k+1)},\cdots,s_n)=\f1\pi\int_{-\?}^\?\f{f(s)}{x_k-s_k}ds_k.
\end{equation}
With this setting, in fact, any finite product of  distinct partial Hilbert transformations may be considered as a Hilbert transformation (also see \cite{car}). Such type of Hilbert transformations corresponds to the Hardy spaces on tubes via the Plemelj formula. Such  Hilbert transformations are of  the tensor product type which are not of interest in this study.

In the present paper we study the Hilbert transformation in the formulation of Clifford algebra. Such formulation is equivalent with the so called conjugate harmonic system studied in \cite{stein} and \cite{sw}.

We note that besides $\RR^n,$ on Lipschitz curves and surfaces, including Lipschitz perturbations of the circle and the real spheres, there exist in depth the $H^\infty$-functional calculi of the respective Dirac type differential operators, or, equivalently, the singular integral operator algebras of monogenic kernels, generalizing the Hilbert transformations in the respective zero-curvature or constant-curvature contexts (\cite{qm}, \cite{q-real ball}, \cite{cq}, \cite{qy}, and references therein).

Now in the Clifford algebra framework  we recall the definition and main properties of the Hilbert transformation representing the conjugate harmonic system. We will be using the Clifford basis elements $e_j, j=1,...,n,$ satisfying $e_j^2=-1, e_je_i=-e_ie_j, i,j=1,...,n, i\ne j.$

\def\ux{\underline{x}}
\def\uy{\underline{y}}
\bT (The Plemelj Theorem) \label{Plemelj theorem} \cite{stein, sw, qm} Let $f$ be a scalar-valued function defined on ${\RR}^n$ in the Clifford algebra setting, $f\in L^p(\RR^n), 1\leq p<\infty.$ Formulate the Cauchy integral of $f$:
\[ Cf(x_0+\ux)=\frac{1}{\omega_{n}}\int_{\RR^n}E(\uy-(x_0+\ux))n(\uy)f(\uy)d\uy, \quad x_0>0, \ \ux\in \RR^n,
\]
where $E(\uy-(x_0+\ux))=\frac{\overline{\uy-(x_0+\ux)}}{|\uy-(x_0+\ux)|^n}$ is the Cauchy kernel in $\RR^{n+1},$ $\omega_{n}$ is the area of the $n$-dimensional unit sphere, $n(\uy)$ is the outer normal at $\uy$ for the upper-half space $x_0>0,$ thus $n(\uy)=-1,$ for any Clifford para-vector $y=y_0+\uy,$ where $y_0\in \RR, \uy\in \RR^n, \overline{y}=\overline{y_0+\uy}=y_0-\uy$ is defined as the Clifford conjugate of $y=y_0+\uy,$ and $d\uy$ is the Lebesgue area measure on $\RR^n.$ When taking limit $x_0\to 0+,$ we have the Plemelj formula
\begin{eqnarray}\label{Plemelj1}
\lim_{x_0\to 0+}Cf(x_0+\ux)=\frac{1}{2}f(\ux)+\frac{1}{2}Hf(\ux), \quad {\rm a.e.},
\end{eqnarray}
where $Hf$ is the Hilbert transform of $f,$ given by
\begin{eqnarray}\label{Hilbert1}
Hf=\sum_{j=1}^n e_jR_jf,
\end{eqnarray}
where $R_j$ is the $j$-th Riesz transformation, being the singular integral operator given by
\begin{equation}\label{Riesz1}
 \nonumber   R_j(f)(\ux)=\lim_{\varepsilon\rt0}\frac{1}{\omega_{n}}\int_{|\uy|>\varepsilon}\f{y_j}{|y|^{n+1}}f(\ux-\uy)d\uy,\
    j=1,\cdots,n.
\end{equation}
\eT

For this result we refer the reader  to \cite{stein}, \cite{sw}, \cite{bg}.

We note that a Clifford-valued smooth function $F$ is said to be left-monogenic, or simply monogenic, if $DF=0,$ where $D=\sum_{k=0}^n \partial_{x_k}e_k,$ with $e_0=1$ (\cite{bds}). We call $x\in \RR^n$ as vector, and in the case write $x$ as $\ux,$ and call $x_0+\ux, x_0\in \RR, \ux\in \RR^n,$ as para-vector. The limit process of the Plemelj Theorem amounts to making the Hardy space projection. By applying twice the Plemelj Theorem we obtain
\[ \left[\frac{1}{2}(I+H)\right]^2=\frac{1}{2}(I+H),\]
 where the last relation implies $H^2=I.$ Taking the limit $x_0\to 0-$ corresponds to making the projection onto the Hardy space in the lower half of the space ${\RR}^{n+1}.$
 \def\uxi{\underline{\xi}}
To state the following theorem we need to introduce the Fourier multipliers $\chi_\pm$  of the Hardy space projection operators (\cite{qm}), namely,
\begin{eqnarray}\label{projection1}
\chi_\pm (\uxi)=\frac{1}{2}\left(1\pm i\frac{\uxi}{|\uxi|}\right).\end{eqnarray}
They satisfy the characteristic properties
\begin{eqnarray}\label{projection2}
\chi_\pm^2=\chi_\pm,\quad \chi_\pm \chi_\mp=0,\quad \chi_++\chi_-=1.\end{eqnarray}

We further have

 \bT \label{equivalence}\cite{stein, sw, qm}  Let $F\in L^2(\RR^n)$ be para-vector-valued. Then the following assertions are equivalent.
\begin{enumerate}
  \item $F$ is the non-tangential boundary limit of a para-vector-valued Hardy $H^2$-function in the upper-half (or the lower-half) space.
  \item  $F$ has the form $F=f+Hf$\ ($F=f-Hf$)\ for some scalar-valued function $f\in L^2(\RR^n).$
  \item $HF=F$ \ (HF=-F).
  \item $\mathfrak{F}(F)=\chi_+\mathfrak{F}(F)$ \ ($\mathfrak{F}(F)=\chi_-\mathfrak{F}(F)$),  where $\mathfrak{F}(F)$ denotes the Fourier transform of $F.$
      \item $F$ can be monogenically extended to the upper- (lower-) half space, still denoted by $F,$ and
   \begin{eqnarray}\label{boundary Cauchy}
   F=CF,\end{eqnarray}
   where $CF$ is the Cauchy integral with boundary data $F.$
   \item $F\ast P_{x_0}$ is a monogenic function in the upper- (lower-) half space, where $P_{x_0}$ is the Poisson kernel at $\pm x_0>0.$
\end{enumerate}
\eT

  We note that in the $L^2(\RR)$ case (that is $n=1$ in the above theorems) with $e_1=-i,$ we have
  $\chi_+=\chi_{(0,\infty)}, \chi_-=\chi_{(-\infty,0)},$ the latter two being, respectively, the characteristic functions of the sets $(0,\infty)$
  and $(-\infty,0),$ and  $H=-i H^0.$ With this relation it is easy to see that the known properties of $H^0$ in relation to the Hardy spaces of the
  upper- and lower-half plane are particular cases with those for general $n$-dimensional spaces described in Theorem
  \ref{equivalence} (\cite{carnett}, \cite{QXYYY}, and therein).

In this paper we prove that the Hilbert transformation defined through boundary limits of monogenic functions may be characterized by
symmetry properties of the operator under motion groups of the space.  This study would exhibit the intrinsic relation between the analyticity
 and the transformation group of the space.  To be more precise, we try to look for the symmetry properties of the underlining operators through the induced
  representations of their space transformation groups similar to ax+b. We discover that the Hilbert transformation $H$ in $\RR^n$ may be characterized by certain symmetry properties in terms of the representations of the motion group of the space. Although $H$ does not commute with rotation, it possesses, and is characterized by, the symmetry of the Spin group, and is invariant under the action of the $r\mathrm{Spin}(n)+\RR^n$ group. In \cite{sw0}, E. M. Stein and G. Weiss point out
    that it is not $\mathrm{SO}(n)$ but $\mathrm{Spin}$ that gives rise to the symmetry of the Dirac equation.  The past studies have revealed certain inter-relationship between the Hilbert transformation, Dirac operator and the Spin group.

 During the past decades there appeared a large number of results concerning induced representation in Clifford analysis (see, for instance, \cite{abb,kq, eh, dm, gm}). In \cite{dm}, V.K.~Dobrev and P.~Moylan, in particular, study the induced representations and invariant integral operators for $SU (2, 2),$ as well as their applications in the coherent state theory. In pursuing our goal we were benefited from the mentioned studies. We give explicit spinor representations of the $r\mathrm{Spin}(n)+\RR^n$ group, by which we deduce a number of symmetry properties of  Hilbert transformations in $\RR^n$. We observe that the theory of the induced representations of the $r\mathrm{Spin}(n)+\RR^n$ group varies greatly along with the dimension number $n$ of the space. In the present paper  we restrict ourselves to $n=2,3$.

For the one-dimensional case, the representation $\pi$ can be
decomposed into two irreducible representations $\pi^\pm$, induced respectively by
the two orbits. They are the only irreducible representations of the trivial group,
over the upper Hardy space and the lower Hardy space.
In the three dimensional case, however,  the two irreducible
representations $\pi^\pm$ are induced  respectively from the two
irreducible representations of the  group $\mathrm{Spin}(n-1)$ on one orbit. In
both the one and three dimensional cases, for $f\in H^+$ we know that
$g(z)=f(\overline{z}),\Re{z}<0,$ is not analytic, i.e., $g(x)\notin
H^-$. The mapping $g(z)\leftrightarrow\overline{f(\overline{z})}$ is
one to one from $H^+$ to $H^-$. But it is not linear since
$\overline{\lambda f(\overline{z})+\mu g(\overline{z})}=\overline{
\lambda}\overline { f(\overline{z})}+\overline{ \mu} \overline {
g(\overline{z})}$. In fact $\pi^+$ is not equivalent to $\pi^-$.It shows that the conjugation is not an isomorphism between the upper and lower Hardy spaces in the case of $n=1,3$

 On the other hand, in the case of $ n = 2$, the conjugation is exactly an isometric isomorphism between the upper and lower Hardy spaces. Denote by
$\rho(f)(\underline{y})=e_2e_1f(-\overline{\underline{y}})$. It can
be easily checked that $\rho$ is an isometric isomorphism from
$H^\pm$ to $H^\mp$, respectively. Moreover we have
$\pi^+\rho =\rho\pi.$
So, $\pi^+$ is isomorphic to $\pi^-$.

The rest of this paper is organized as follows: In \S 2, we introduce the $r\mathrm{Spin}(n)+\RR^n$ group, and briefly describe the symmetry of the Riesz and the Hilbert transformations. In \S 3, we introduce the basic induced representations of the $r\mathrm{Spin}(n)+\RR^n$ group, emphasizing the cases $n=2$ and $3$. In \S 4, we propose the other equivalent representations of   $r\mathrm{Spin}(3)+\RR^3$  in the setting of $\HH$ and $\mathrm{Cl}_3$ separately, and one for  $r\mathrm{Spin}(2)+\RR^2$. \S 5, we characterizes the Hilbert transformation by the representation in the setting of $\HH$. In \S 6, we continue to characterize the Hilbert transformations in the settings of $\mathrm{Cl}_3$ and $\mathrm{Cl}_2.$  \S 7 states that if a group acting on $\RR^n$ naturally commutes with the Hilbert transformation, then it is essentially the  $r\mathrm{Spin}(n)+\RR^n$ group.

\section{The symmetry of Hilbert transformation in $\RR^n$}

It is pointed out in \S 1 that on $\RR$ the Hilbert
transformation commutes with the translation and dilation, and then
it is invariant under the natural action of the ax+b group over
$L^2(\RR)$. In this sense we say that the Hilbert transformation has
the symmetry of ax+b group.
 We studied the Hilbert transformation on the real line and that on the unit circle relation to their symmetry properties in \cite{dlq}.
In the present paper we study the Hilbert transformations in $\RR^n$ relation to their motion group symmetry properties.

There already exist some analogous results concerning symmetry of the Riesz transformations $R_j$ in $\RR^n$,\ $ j=1,2,\cdots,n.$

 \bT \label{riesz}\cite{stein}  Let $T=(T_1,T_2,\cdots,T_n)$
be a n-tuple of bounded transformations on $L^2(\RR^n)$. Suppose
\begin{enumerate}
  \item Each $T_j$ commutes with the translation of
$\RR^n$,
  \item  Each $T_j$ commutes with the dilations of
$\RR^n$,
  \item For every rotation $\rho=(\rho_{jk})$ of $\RR^n$,
$\rho T_j\rho^{-1}=\sum_k\rho_{jk}T_kf$.
\end{enumerate}
\indent Then $T$ is a constant multiple of the $n$-tuple of the Riesz
transformations, i.e. there exists a constant $c$ such  that $T_j=cR_j$,
$j=1,\cdots,n$.\eT

Theorem 2.1 characterizes the Riesz system in terms of their symmetry properties. But it does not precisely correspond to the type of results of the one dimensional case in term of the representations of the motion groups.

It turns out that in the Clifford algebra language one can spell out the right symmetry properties characterizing the Hilbert transformation in $\RR^n$.

Let $e=\{e_1,e_2,\cdots,e_n\}$ be an orthogonal basis for the real
Clifford algebra $\mathrm{CL}_n$
as defined before. We identify $\RR^n$ with the subspace of
$\mathrm{CL}_n$ , the real linear subspace
$\mathrm{span}\{e_1,e_2,\cdots,e_n\}$ of $\mathrm{Cl}_n$:
$x=(x_1,x_2,\cdots,x_n)\rightarrow
\underline{x}=x_1e_1+x_2e_2+\cdots+x_ne_n$. The notation
$\underline{x}$ is sometimes denoted as just $x.$

  Denote by
$L^2(\RR^n,\mathrm{Cl}_n)$ the set of the square integrable
Clifford-valued functions. Then by (\ref{hilbert1}), the Hilbert
transformation $H$ over $L^2(\RR^n,\mathrm{Cl}_n)$ is defined by
\begin{equation}\label{}
\nonumber H(f)(\ul{x})=\lim_{\varepsilon\rt0}\f1{\omega_n}\int_{|\ul{y}|>\varepsilon}\f{\ul{y}}{|\ul{y}|^{n+1}}f(\ul{x}-\ul{y})d\ul{y},\
\end{equation}
for $f\in L^2(\RR^n,\mathrm{Cl}_n)$. As a convolution type singular integral
operator,  $H$ may be alternatively written as
\begin{equation}\label{hilbert3}
H(f)(\ul{x})=c_n\int_{\RR^n}\f{\ul{x}-\ul{y}}{|\ul{x}-\ul{y}|^{n+1}}f(\ul{y})d\ul{y},
\end{equation}
where we suppress the abbreviation p.v. for principal value of the singular integral.

The Hilbert transformation $H$ defined by (\ref{hilbert3})
is invariant under both dilation and translation. It is unexpected that
it can not commute with the natural rotation action.  For
example, let $A\in \mathrm{SO}(n)$ be an  orthogonal matrix with
determent $1$. Then it naturally acts over
$L^2(\RR^n,\mathrm{Cl}_n)$ as $\mathcal{A}f(x)=f(Ax)$. Notice that
\begin{eqnarray}\label{orth}
H(\mathcal{A}f)(x)&=&c_n\int_{\RR^n}\f{\ul{x}-\ul{y}}{|\ux-\uy|^{n+1}}f(A\uy)d\uy\\
 \nonumber&=&c_n\int_{\RR^n}\f{\ul{x}-\ul{A^{-1}\eta}}{|\ux-A^{-1}\eta|^{n+1}}f(\underline{\eta})d(A^{-1}\underline{\eta}),
 \end{eqnarray}
 where $\underline{\eta}=A\uy$.
Although both the Lebesgue measure and norms of vectors are invariant
under rotation, it can be checked that (\ref{orth}) does not equal
to $\mathcal{A}(Hf)$ for  $f\in L^2(\RR^n,\mathrm{Cl}_n)$ in general. Thus the Hilbert transformation $H$ does not possess the symmetry of $\mathrm{SO}(n)$.

It also worthies mentioning that in the space $\RR^3$ the operator $Hf$ is not invariant under the action of reflection through a single plane since both the Cauchy and the Hilbert integral are of the second class (\cite{CZ}).
\vspace{2mm}

We need to involve an equivalent definition of $\mathrm{SO}(n)$. Consider the group
$$
\mathrm{Spin}(n)=\{\ul{u}\ul{v}:\ul{u},\ul{v}\in \RR^{n}\ \mathrm{
and }\ |\ul{u}|=|\ul{v}|=1\}.
$$
The spin group is the double covering of the special orthogonal group $\mathrm{SO}(n)$. Now denote
by $r\mathrm{Spin}(n)+\RR^n$ the triplet $(r,s,b)$ with $r>0$, $s\in
\mathrm{Spin}(n)$ and $b\in\RR^n$.  Assign the action of an element $(r,s,b)$ of $r\mathrm{Spin}(n)+\RR^n$
over $\RR^n$ as
\begin{equation}\label{rspin}
    (r,s,b)\ul{x}=r(s\ul{x}s^{-1})+b.
\end{equation}
Denote $\chi((1,s,0))\ul{x}=s\ul{x}s^{-1}$. By \cite{dss}, for
any $\ul{x}\in\RR^n$ and $s\in \mathrm{Spin}(n)$, the Clifford number $s\ul{x}s^{-1}$
 belongs to $\RR^n$. So (\ref{rspin}) is well defined as a linear map from $\RR^n$ to itself. For
$(r,s,b),(r',s',b')\in r\mathrm{Spin}(n)+\RR^n$ their composition
also belongs to $r\mathrm{Spin}(n)+\RR^n$ since
\begin{equation}
\nonumber
(r,s,b)((r',s',b')\ul{x})=(r,s,b)(r's'\ul{x}s'^{-1}+b')=rs(r's'\ul{x}s'^{-1}+b')s^{-1}+b,\
\ul{x}\in\RR^n.
\end{equation}
Hence $r\mathrm{Spin}(n)+\RR^n$ is the group by the multiplication law
\begin{equation}\label{grouplaw}
(r,s,b)(r',s',b')=(rr',ss',rsb's^{-1}+b).
\end{equation}
We denote by $\mathfrak{G}_n$ this group and call it the
$r\mathrm{Spin}(n)+\RR^n$ group. It is a subgroup  of the universe
covering of  the affine group of $\RR^n$.

The natural unitary representation $\pi$ of $\mathfrak{G}_n$ over $
L^2(\RR^n,\mathrm{Cl}_n)$  is defined by
\begin{equation}\label{spinrepn}
    (\pi(r,s,b))f(\ul{x})=(\f1r)^{\f n2}sf(\f1r s^{-1}(\ul{x}-b)s),\
    (r,s,b)\in\mathfrak{G}_n,f\in L_2(\RR^n,\mathrm{Cl}_n),
\end{equation}
where $\pi(r,s,b)$ is sometimes written as $\pi_{rsb}$.
\vspace{3mm}

\noindent The Hilbert transformation can naturally be linearly extended to
act on $L^2(\RR^n,\mathrm{Cl}_n)$ because $
L^2(\RR^n,\mathrm{Cl}_n)$ is the direct sum of $2^{n}$ spaces
isomorphic to $L^2(\RR^n)$. Then we can state the symmetry of the
Hilbert transformation in $\RR^n$.

\bP\label{p1} The Hilbert transformation commutes with the natural
representation of  $\mathfrak{G}_n$.\eP

 \bpp
It can be checked by a direct computation. \epp

We need the so called Dixmier and Schur's Lemmas.

\bT\cite{kn} {\bf(Dixmier's Lemma)}\label{schur} Let $\mathcal{H}$ be a complex Hilbert space. If a
family of transformations $\mathcal{A}$ acts irreducibly on
$\mathcal{H}$, then any linear transformation $C$ that commutes
with all $T\in \mathcal{A}$ is of the form $C=\lm\mathcal{I}$, where $\lm$ is
a complex number and $\mathcal{I}$ is the identity transformation.\eT

Let $\s$ be a representation of the group $G$ over the Hilbert space
$\mathcal{H}$. It is obvious that $\s$ is irreducible
  if and only if $\{\s(g):g\in G\}$ acts irreducibly on
$\mathcal{H}$. As a consequence of the last theorem the following result holds.

\bT {\bf(Schur's Lemma)}\label{t1} Suppose that $\pi'$  is an
irreducible sub-representation of $\pi$ over $\mathcal{H}'$, a
closed subspace of $L^2(\RR^n, \mathrm{Cl}_n)$. Then over $\mathcal{H}'$ the Hilbert transformation H acts as $H=\lm\mathcal{I},$ where $\lambda$ is some complex number. \eT

In the rest of the paper we try to find the invariant subspaces of
$\pi$, and then characterize the operator $H$.

\section{induced representation of $r\mathrm{Spin}(n)+\RR^n$ group}
In this section we give the { induced} representations of the
$r\mathrm{Spin}(n)+\RR^n$ group from $\mathrm{Spin}(n-1)$'s representations. In \cite{sa} S. Sahi and E.M. Stein obtained its explicit representation on the Hilbert matrix space.   In this paper we obtain their representations
induced by the basic spinor representations of
$\mathrm{Spin}(n)$.

 \subsection{Examples of $\mathrm{Spin}(n)$} We first recall  some fundamental properties of $\mathrm{Spin}(n)$ and its representation (See \rm \cite{dss} for more details).
 \vspace{2mm}

 \bE \label{spin1} For $n=1, \RR^{0,1}=\RR e$ with $e^2=-1$, then $\mathrm{Spin}(1)=\{1,-1\}\cong\ZZ_2.$ \eE

 \bE\label{spin2} $\mathrm{Spin}(2)\cong U(1)\cong  \mathrm{SO}(2).$
 It is obvious that  $  \mathrm{Spin}(2)$ is a multiplicative subgroup of Clifford algebra $\mathrm{Cl}_2$, each of whose  elements $s$ can be written as
\begin{equation}\label{sp2}
  s=\cos\th+\sin\th e_{1}e_2,
\end{equation}
where $\th$ is just the angle of the rotation represented by $s$. \eE

\bE\label{spin3}
For $n=3$ let $s\in\mathrm{Spin}(3)$ and  $\ul{\om},\ul{\nu}\in S^2$ such
 that $s=\ul{\om}\ul{\nu}$, and $\th\in[0,\pi]$ is the angle between $\ul{\om}$ and $\ul{\nu}$. Then we have that
\begin{equation}\label{spin-act}
  \chi(s)(\ul{x})=s\ul{x}s^{-1},\ x\in \RR^3,
\end{equation}
is the rotation determined by $s$  around the axis
$\ul{\om}\wedge\ul{\nu}e_{123}$ and through the angle $2\th$ in
$\RR^3$. \vspace{2mm}

 $\mathrm{Spin}(2)$ can be viewed as identical with the subgroup of $\mathrm{Spin}(3)$ in which every $s\in\mathrm{Spin}(2)$ is considered as a
 rotation in $\RR^3$ fixing $e_3$. We have $\mathrm{Spin}(3)\cong S^3$ and
\begin{equation}\label{spin3-spin2}
  \mathrm{Spin}(3)/\mathrm{Spin}(2)\cong S^3/S^1\cong S^2.
\end{equation}

\eE

We give a link between any two elements of $\mathrm{Spin}(3)$ by
$\mathrm{Spin}(2)$ (\cite{gw}). For $s,t\in \mathrm{Spin}(3)$, through
identifying $t$ and $s$ with, respectively,
$$
\begin{array}{l}
  \left(
     \begin{array}{cc}
       z' & w' \\
       -\ol{w}' & \ol{z}' \\
     \end{array}
   \right)
    \ \text{and}\ \left(
                     \begin{array}{cc}
                       z & w \\
                        -\ol{w}& \ol{z} \\
                     \end{array}
                   \right),
\end{array}
$$
    the requirement $t=su$ for some $u=\cos\f\theta2+\sin\f\theta2 e_{12}\in\mathrm{Spin}(2)$ then gives rise to

 \begin{equation}\label{matrix2}
  \begin{array}{l}
\left(
     \begin{array}{cc}
       z' & w' \\
       -\ol{w}' & \ol{z}' \\
     \end{array}
   \right)=  \left(
                     \begin{array}{cc}
                       z & w \\
                        -\ol{w}& \ol{z} \\
                     \end{array}
                   \right)\left(
                            \begin{array}{cc}
                              e^{i\f\th2} &  \\
                               & e^{-i\f\th2} \\
                            \end{array}
                          \right),
\end{array}
\end{equation}
 or, alternatively (\cite{dss}),
\begin{eqnarray}
 \nonumber  z'&=&ze^{i\f\th2}=\cos\a e^{i(\varphi+\f\th2)},  \\
\nonumber  w' &=& we^{-i\f\th2}=\sin\a e^{-i(\psi+\f\th2)}.
\end{eqnarray}

\subsection{Representation of $\mathrm{Spin}(n)$}
 The representations of the groups $\mathrm{Spin}(n)$ may
be derived from the representations of complex Clifford algebras on the
\bf spinor spaces. \rm  Denote by $\CC^{2m}$ the complexifications
of $\RR^{0,2m}$, whose associated complex Clifford algebra is
$\CC_{2m}$. Let $\{e_{2j-1},e_{2j}:j=1,\cdots,m\}$ be an orthonormal
basis of $\CC^{2m}$. Define, for each $j=1,\cdots,m$, the elements
$f_j$ and $\ol{f}$ by, alternatively, $ f_j= \f12(e_{2j-1}+ie_{2j})$ and $
\ol{f}_j=-\f12(e_{2j-1}-ie_{2j}). $ Let
$W=\mathrm{Span}_\CC\{f_1,\cdots,f_m\}$ and
$\ol{W}=\mathrm{Span}_\CC\{\ol{f}_1,\cdots,\ol{f}_m\}$. Denote by
$\wedge W$ and $\wedge \ol{W}$ the Grassmann algebras determined by
$W$ and $\ol{W}$, respectively. Set $
  I=\prod_{j=1}^m\ol{f}_jf_j.$ Then
 $(\wedge W)I=\CC_{2m}I$ is a minimal
left ideal in $\CC_{2m}$.

 Let $S_{2m}=(\wedge W)I$ and $S^+_{2m}=(\wedge^{\text{even}} W)I$,
$S^-_{2m}=(\wedge^{\text{odd}} W)I$.
 \bE \label{spinor2} The case
of $m=1$: $f_1=\f12(e_1+ie_2),\ol{f}_1=-\f12(e_1-ie_2)$  and
$I=\ol{f}_1f_1$, and $S_2=(\CC\oplus\CC f_1)I$,  $S^+_2=\CC I$ and
$S_2^-=\CC f_1I$. \eE
The left multiplication of $\mathrm{Spin}(2m)$
induces an irreducible representation on $S^+_{2m}$.

\bT\cite{dss}\label{representation} Let the actions of the group
$\mathrm{Spin}(2m)$  be the left multiplication in $\CC_{2m}$. We
have two irreducible complex representations of
$\mathrm{Spin}(2m)$ on, respectively, $S^+_{2m}$ and $S^-_{2m}$.\eT
\bT\cite{dss}\label{representation1} Let the actions of the group
$\mathrm{Spin}(2m-1)$  be the left multiplication in $\CC_{2m}$. We
have two equivalent irreducible complex representations of
$\mathrm{Spin}(2m-1)$ on, respectively, $S^+_{2m}$ and $S^-_{2m}$. \eT

\subsection{Induced Representation of $r\mathrm{Spin}(n)+\RR^n$ }
 Now we discuss  the structure of the $r\mathrm{Spin}(n)+\RR^n$ group
$\mathfrak{G}_n$. $\mathrm{Spin}(n)$ is  underlain by a
compact manifold $S$ of dimension $\f{n(n-1)}2$. Denote
$M=\{(r,s,0):r>0,s\in\mathrm{Spin}(n)\}$. Then $M$ is a closed
subgroup of $\mathfrak{G}_n$, whose underlying manifold is
$(0,\?)\times S$.  $M$ is the direct product of $\mathrm{Spin}(n)$ and
  $\RR^+$ with the group law as standard multiplication, i.e.,
  $M\cong\RR^+\oplus\mathrm{Spin}(n)$. Again denote  $N=\{(1,1,b):b\in \RR^n\}$ which is  a closed subgroup.
Let $(r,s,b)\in M$ and $(1,1,b')\in N$. Then
$(r,s,b)^{-1}=(\f1r,s^{-1},-\f1r s^{-1}bs)$ and
$$
(r,s,b)^{-1}(1,1,b')(r,s,b)=(\f1r,s^{-1},-\f1r
s^{-1}bs)(r,s,b'+b)=(1,1,\f1r s^{-1}b's)\in N,
$$
which amounts to the fact  that $N$ is a normal subgroup of $\mathfrak{G}_n$.
   It is easy to check that the map $(n,h)\rt nh$ from $N\times M$ to $\mathfrak{G}_n$
is a group isomorphism. Hence $\mathfrak{G}_n$ is the \bf semi-direct
product \rm of $N$ and $M$, \cite{folland}, denoted by
$\mathfrak{G}_n=N\ltimes
H=\RR^n\ltimes(\RR^+\oplus\mathrm{Spin}(n))$. In this paper we regard $\mathfrak{G}_n$ as the $r\mathrm{Spin}(n)+\RR^n$ group so as to remind the similarity with the \textbf{ax+b group}.

\vspace{2mm}
In the sequel  we equally use all the three notations, namely,
$r\mathrm{Spin}(n)+\RR^n$,
$\RR^n\ltimes(\RR^+\oplus\mathrm{Spin}(n))$ and $\mathfrak{G}_n$. The notation
$\mathfrak{G}_n$  is sometimes abbreviated as $\mathfrak{G}$. Since we can not find a commonly accepted terminology for the group $r\mathrm{Spin}(n)+\RR^n$ in the reachable references, we sometimes also call it the ${\bf ax+b \rm \ group}.$

\vspace{2mm}

Now we introduce the Mackey machinery  \cite{ka} for the induced representation of  $\mathfrak{G}_n$.

As a consequence of the fact that the additive group $\RR^n$ is an Abel group, any irreducible unitary representation of $\RR^n$ has dimension one.
 Let $\rho$ be
a unitary representation of $\RR^n$. Then there exists $\xi\in
\RR^n$ such that $\rho(1,1,x)=e^{i2\pi x\cdot\xi}$. Thus we may
identify the dual group of $N$ with  $\RR^n$. Let $\mathfrak{G}$ act
on $N$ by conjugation. It induces an action of $\mathfrak{G}$
on the dual group $\hat{N}$, $(x,\nu)\rt x\circ\nu$, defined by
\begin{equation}\label{stab}
    <n,x\cdot\nu>=<x^{-1}nx,\nu>\ (x\in \mathfrak{G},\nu\in\hat{N},n\in
    N),
\end{equation}
where $<x,\nu>$ is the action of $\nu$ on $x$ as a linear functional. Without ambiguity we sometimes
 denote $x\cdot\nu$ as $x\nu$ in the sequel. For each $\nu\in
\hat{N}$, we denote by $\mathfrak{G}_\nu$ the stabilizer of $\nu$,
\begin{equation}\label{}
    \nonumber \mathfrak{G}_\nu=\{x\in \mathfrak{G}:x\nu=\nu\},
\end{equation}
which is a closed subgroup of $\mathfrak{G}$. We denote by
$\mathcal{O}_\nu$ the orbit of $\nu$:
\begin{equation}\label{}
    \nonumber \mathcal{O}_\nu=\{x\nu: x\in \mathfrak{G}\}.
\end{equation}
For $x=(r,s,b)\in \mathfrak{G}$, recalling that $\chi(s)$ is a rotation, by (\ref{stab}), we have
\begin{eqnarray}
  \nonumber<n,x\nu>&=&<(\f1r,s^{-1},-s^{-1}bs)(1,1,n)(r,s,b),\nu>\\
  \nonumber&=&<(1,1,\f1r s^{-1}ns),\nu>\\
 \label{stab1} &=&<\f1r s^{-1}ns,\nu>=<n,\f1r s\nu s^{-1}>,\ x\in\mathfrak{G},\nu\in\hat{N},n\in N,
\end{eqnarray}
where $(1,1,n)$ is made to be identified with $n$. So
\begin{equation}\label{orbit}
x\nu=\f1r s\nu s^{-1}.
\end{equation}
Identify $\nu$ with $(1,1,\nu)$  and let
$\nu=(\nu_1,\nu_2,\cdots,\nu_n)$. Then $x\nu=\nu$ if and only if
$r=1,$ where $\chi(x)$ stand for rotations around the axis $\nu$. Because
the subgroup of $\mathrm{Spin}(n)$ about the fixed axis is
isomorphic to $\mathrm{Spin}(n-1)$, we obtain the relation
\begin{equation}\label{stab2}
  \mathfrak{G}_\nu=\RR^n\ltimes\mathrm{Spin}(n-1).
\end{equation}
\vskip 3mm
Then we can completely classify the irreducible
representations of $\mathfrak{G}_\nu$ in terms of the irreducible
representations of $\RR^n$ (i.e., the characters $\nu\in \hat{N}$) and
the irreducible representations of their \bf little groups \rm $M_\nu$ associated with $\nu$:
\begin{equation}\label{}
    \nonumber M_\nu=\mathfrak{G}_\nu\cap M.
\end{equation}
By (\ref{stab2}) we obtain that $M_\nu=\mathrm{Spin}(n-1)$. For
example, $M_\nu=\ZZ_2$ for $n=2$.

 \bL \cite{folland}\label{gv}
Let $G_\nu$ be a locally compact group and $N$ its normal Abel subgroup, and $G_\nu=N\ltimes M_\nu $.  Suppose that $\rho$ is an irreducible representation of $M_\nu$. Then we
have an irreducible representation of $G_\nu$, denoted by
$\nu_\rho$, satisfying
\begin{equation}\label{}
    \nu_\rho(nh)=<n,\nu>\rho(h),\forall n\in N,h\in M_\nu.
\end{equation}
On the contrary  every irreducible representation of $G_\nu$ can be written into such a form.
\eL
For $n\geq3$, $\mathrm{Spin}(n-1)$ acts transitively on the unit
sphere $S^{n-1}$. So by (\ref{orbit}), $\mathfrak{G}$ acts
transitively on $\hat{N}$. Then there exist only two orbits:
$\mathcal{O}_0$, $\mathcal{O}_{e_n}=\mathcal{O}_{(0,0,\cdots,1)}$.

By Example \ref{spin1}, $\mathrm{Spin}(1)=\ZZ_2$, that does not act transitively
on the unit circle. But $\RR^2\ltimes(\RR^+\oplus\mathrm{Spin}(2))$
does on
 $\hat{\RR^2}$. There  exist  two orbits:  $\mathcal{O}_0$, $\mathcal{O}_{e_2}=\mathcal{O}_{(0,1)}$.
 Being contrary to the case $n\geq3$,  it is obvious that $\mathrm{Spin}(1)$  contributes little to the transitivity of $\RR^2\ltimes(\RR^+\oplus\mathrm{Spin}(2))$, which underlines the complexity of the representation.

According to the definition given in \cite{folland},  $\mathfrak{G}$ can
  be said to acts \bf regularly \rm on $\hat{N},$
because the Borel set $\{0,e_n\}$ only intersects $\mathcal{O}_0$
and $\mathcal{O}_{e_n}$ at $0$ and $e_n,$ respectively. \vspace{2mm}

\bT\cite{folland} \label{main1} Suppose $G=N\ltimes M$, where $N$ is
Ableian and $G$ acts regularly on $\hat{N}$. If $\nu\in\hat{N}$
 and $\rho$ is an irreducible representation of $M_\nu$, then the induced representation
 $\mathrm{ind}^G_{G_\nu}(\nu\rho)$ is an irreducible representation
 of $G$, and every irreducible representation of $G$ is equivalent
 to one of those forms. Moreover, $\mathrm{ind}^G_{G_\nu}(\nu\rho)$ and
 $\mathrm{ind}^G_{G_{\nu'}}(v'\rho')$ are equivalent if and only if
 $\nu$ and $\nu'$ belong to the same orbit. That is, if for some $x$ there holds $\nu'=x\nu$, then $h\rt
 \rho(h)$ and $h\rt\rho'(x^{-1}hx)$ are equivalent representations
 of $M_\nu$.\eT
\subsection{Cases $n=2,3$}
In  this section we study  the induced
representations of $\mathfrak{G}$ for $n=2,3$ under the process suggested
by Theorem \ref{main1}.

 By (\ref{stab2}) we have $
\mathfrak{G}_\nu=\RR^n\ltimes\mathrm{Spin}(n-1)$. By Example
\ref{spinor2} there exist representations of $\mathrm{Spin}(2)$ over $S^+_2$ and $S^-_2$.

\bE \label{gnu2}When $n=2$,
$\mathfrak{G}_\nu=\mathrm{Spin}(1)+\RR^2=\RR^2\ltimes\ZZ^2$.
It is obvious that representations of $\mathfrak{G}_\nu$ should
be of the forms
\begin{equation}\label{rho0}
    \rho(s,b)x=e^{i2\pi t\cdot {b}}sx,\ (s,b)\in\mathfrak{G}_\nu, x\in S^+_2
\end{equation}
or
\begin{equation}\label{rhoo}
    \rho(s,b)x=e^{i2\pi t\cdot {b}}sx,\ (s,b)\in\mathfrak{G}_\nu,x\in S^-_2,
\end{equation}
where $t\in\RR^2$.

These two representations are equivalent to each other.
 \eE

\bE \label{gnu3}Denoted by $\rho^+$ and $\rho^-$ the representations
of $\mathrm{Spin}(2)$ over $S^+_2$ and $S^-_2$, respectively. Let
$\nu$ be an irreducible unitary representation of $\RR^3$, i.e.,
there is $t\in\RR^3$ such that $\nu(b)=e^{i2\pi t\cdot {b}}$ for
$b\in \RR^3$.
If $\rho$ is an irreducible representation  of
$\RR^3\ltimes\mathrm{Spin}(2)$, then there exists $t\in\RR^3$  such that
either
\begin{equation}\label{rho1}
\rho(s,b)=e^{i2\pi t\cdot {b}}\rho^+(s)
\end{equation}
 over the Hilbert space $S^+_2$, or
\begin{equation}\label{rho2}
\rho(s,b)=e^{i2\pi t\cdot {b}}\rho^-(s)
\end{equation}
 over the Hilbert space $S^-_2$.

We note that $\rho^\pm(s)$ are respectively the left multiplication by $s$ over
$S^\pm_2.$

 \eE

There  does not always exist an invariant measure for a noncompact
non-Abel group. But $dx=\f1rdr\times d\s(s)\times dm(b),$ for
$x=(r,s,b)\in\mathfrak{G}$, is the left invariant measure for
$\mathfrak{G}=\RR^n\ltimes(\RR^+\oplus\mathrm{Spin}(n)),n=2,3$.

Now we can present the representations of
$\mathfrak{G}=\RR^n\ltimes(\RR^+\oplus\mathrm{Spin}(n))$ induced by
the representations of $\mathrm{Spin}(n-1)$. Notice that here $\mathrm{Spin}(n-1)$,
considered as the collection of the spins fixing the axis $\nu=e_n$, is
a subgroup  of $\mathrm{Spin}(n)$. We note that in the  definition of $\nu(b)$ the
selected vector $t$ should take $t=e_n.$
Recall that $\mathfrak{G}_\nu=\RR^n\ltimes\mathrm{Spin}(n-1)$ is a closed
subgroup of the locally compact group $\mathfrak{G}$. Then $\mathfrak{G}/ \mathfrak{G}_\nu$ is also  locally compact. By the homeomorphism
$\mathfrak{G}/ \mathfrak{G}_\nu\cong \RR^+\times S^n$, it is easy to check that $\mathfrak{G}/ \mathfrak{G}_\nu$ admits an invariant
 measure $\mu$ such that $d\mu(x\mathfrak{G}_\nu)=\f1r drd\Sigma$, where $\Sigma$
  is the Borel measure of the sphere $S^n$.
\vskip3mm
  For $n=3$, by Example \ref{representation}
  we get  that there is an irreducible unitary representation $\rho^+$
of $\mathfrak{M}_\nu=\mathrm{Spin}(2)$ on the Hilbert space $S^+_2$
with the standard inner product. By Example \ref{gnu3} there exists
one and only one  irreducible representation $\rho$ of
$\mathfrak{G}_\nu$. Then we can induce a
representation of $\mathfrak{G}$ by $\rho^+.$ Denote by
$\mathrm{C}(\mathfrak{G},S^+_2)$ the space of continuous Weyl
spinor-valued functions from $\mathfrak{G}$ to $S^+_2$. To introduce the representation space we
begin with introducing the function space
\begin{equation}\label{F0}
   \begin{array}{rl}
    {\mathcal F}_0
    =\{f\in \mathrm{C}(\mathfrak{G},S^+_2):\quad &\text{the support of }\ f\ \text{is compact  and}\\[1mm]
     &f(x\z)=\rho(\z^{-1})f(x),\ \forall x\in \mathfrak{G},\ \z\in \mathfrak{G}_\nu\},
     \end{array}
\end{equation}
where $x=(r,s,b)\in r\mathrm{Spin}(3))+\RR^3$,
$\z=(1,\eta,c)\in\mathrm{Spin}(2)+\RR^3$ and
$\rho(\z^{-1})=\rho^{-1}(\z)=e^{-i2\pi \nu\cdot
{c}}(\rho^+({\eta}))^{-1}$. Recalling that $\rho^+(\eta)$ is just the left multiplication of $\eta$, the equality in (\ref{F0}) may be
rewritten as
\begin{equation}\label{coset}
  f((r,s,b)(1,\eta,c))=e^{-i2\pi {c}\cdot\nu}\eta^{-1}f(r,s,b).
\end{equation}

 By (\ref{F0}) and (\ref{coset}) we can define a function
$\tilde{f}(\tilde{x})$ over $\mathfrak{G}/\mathfrak{G}_\nu$ by
$\tilde{f}(\tilde{x})=f(x)$ for $\tilde{x}=x\mathfrak{G}_\nu\in
\mathfrak{G}/\mathfrak{G}_\nu$. That is, $f$ is essentially  the
function over $\mathfrak{G}/\mathfrak{G}_\nu$. In fact, by Proposition 6.1 in \cite{folland}, for every element $f$ of ${\mathcal
F}_0$ there exists  a Weyl spinor-valued function
$\a_f:\mathfrak{G}_\nu\rt S^+_2$, being continuous  with compact
support such that
\begin{equation}\label{gh}
    f(x)=\int_{\mathfrak{G}_\nu}\rho(\eta)\a_f(x\eta)d\eta,\quad x\in
    \mathfrak{G},
\end{equation}
where $d\eta$ is the invariant measure of $\mathfrak{G}_\nu$.

Denoted by $<\cdot,\cdot>_0$ the inner product of the complexification
of $\RR_{0,n}$. By (\ref{gh}), for $f,g\in{\mathcal F}_0$,
$<f(x),g(x)>_0$ depends only on the $\mathfrak{G}_\nu$-coset of $x$.
Thus $<f(x),g(x)>_0$ can be considered as a complex function on
$\mathfrak{G}/\mathfrak{G}_\nu$. Thus we can define an inner product
on ${\mathcal F}_0$ as
\begin{equation}\label{product}
    <f,g>=\int_{\mathfrak{G}/\mathfrak{G}_\nu}<f(x),g(x)>_0 d\mu(x \mathfrak{G}_\nu).
\end{equation}

 Let ${\mathcal F}$ be the
Hilbert space completion of ${\mathcal F}_0$. Then we  get a unitary
representation $\s$ of $\mathfrak{G}$ on ${\mathcal F}$ by the left
translation $L_x,\ x\in G$,
\begin{equation}\label{left}
    L_x(f)(y)=f(x^{-1}y),\ \forall y\in \mathfrak{G},
\end{equation}
that is, $\s(x)=L_x$. Then $\s$ is called the induced representation
by $\rho$ and denoted by
$\mathbf{Ind}^\mathfrak{G}_{\mathfrak{G}_\nu}(\rho)$.
\vspace{3mm}

\noindent We need parameterizations of  $\mathfrak{G}/\mathfrak{G}_\nu$ and that of the elements of ${\mathcal F}$.
 Since
\begin{equation}\label{}
\nonumber\mathfrak{G}/\mathfrak{G}_\nu=(\RR^3\ltimes(\RR^+\oplus\mathrm{Spin}(3)))/(\mathrm{\RR^3\ltimes Spin}(2)),
\end{equation}
we obtain, by (\ref{spin3-spin2}),
\begin{equation}\label{paramet}
\mathfrak{G}/\mathfrak{G}_\nu\cong \RR^+\times
(\mathrm{Spin}(3)/\mathrm{Spin}(2))\cong \RR^+\times S^2,
\end{equation}
where $\cong$ denotes the isometric homeomorphism.

\vspace{2mm}
\noindent From now on, the symbol $x$ does not stand for  an element of $\mathfrak{G}$ but a vector of $\RR^3$.\\
\vspace{2mm}

\noindent By Example \ref{spinor2} we have
\begin{eqnarray}\label{hernion}
\nonumber S^+_2 &=&  \CC I=\CC(\ol{f}_1 f_1)=\CC(-\f12(e_1-ie_2)\f12(e_1+ie_2))\\
  &=&\CC(1-ie_1e_2).
\end{eqnarray}

An element $f$ of $\mathcal{F}_0$ may be written as
\begin{equation}\label{}
\nonumber  f(r,s,b)=\a(r,s,b)(1-ie_1e_2), (r,s,b)\in\mathfrak{G},
\end{equation}
where $\a(r,s,b)$ is a scalar-valued function. Then $f$ satisfies the ralation
\begin{eqnarray}\label{firsttrans}
 \nonumber  f((r,s,b)(1,\eta,c))&=&  e^{-i2\pi c\cdot {\nu}}\eta^{-1}\a(r,s,b)(1-ie_1e_2),\\
   &\ &(r,s,b)\in\mathfrak{G},(1,\eta,c)\in \mathfrak{G}_\nu.
\end{eqnarray}

By (\ref{spin3-spin2}), there exists an isometrical homeomorphism
between   $\mathrm{Spin}(3)/\mathrm{Spin}(2)$ and the unit sphere $S^2$ by the correspondence $s\rightarrow \om =s\nu s^{-1}$.
For $\om\in S^2$, choose $s_\om\in \mathfrak{G}$ such that  $\tilde{s}_\om=s_\om
\mathfrak{G}_\nu$ is its corresponding element in
$\mathfrak{G}/\mathfrak{G}_\nu$, i.e., ${\om}=s_\om e_3s^{-1}_\om$. Recall that $s$ admits a matrix representation of (\ref{sp2}), where $z=\cos\a e^{i\varphi}, w=\sin\a
e^{-i\psi},0\<\a\<\f\pi2,\varphi,\psi\in\RR$. By(\ref{matrix2}) we may assume that $s_\om$ has a matrix representation as
\begin{equation}\label{matrix3}
S_\om=\left(
                                \begin{array}{cc}
                                  \cos\beta  & \sin\beta e^{-i(\psi-\varphi)} \\
                                  {\sin\beta e^{i(\psi-\varphi)}} & {\cos\beta } \\
                                \end{array}
                              \right),
\end{equation}
i.e., we choose $\theta=2\varphi$. Let $\left(
                                   \begin{array}{cc}
                                     e^{i\f\th2} & 0\\
                                     0 & e^{-i\f\th2} \\
                                   \end{array}
                                 \right)$ be  the matrix representation of $\eta_{\om s}\in\mathrm{Spin}(2)$. Then  $s=s_\om \eta_{\om s}$.

Recall that $(r,s,b)(1,\om,c)=(r,s\om,rscs^{-1}+b)$. If we take
$b=0,s=s_\om,\eta=\eta_{\om s}$ and
write $\f1r s^{-1}bs$ in place of $c$ in (\ref{firsttrans}), we have
\begin{eqnarray}\label{firsttrans'}
 f(r,s,b) &=&f((r,s_\om,0)(1,\eta_{\om s},\f1r s^{-1}bs)\\
 \nonumber&=& e^{-i2\pi \f1r {s^{-1}bs}\cdot{\nu}}\eta_{\om s}^{-1}\a(r,s_\om,0)(1-ie_1e_2).
\end{eqnarray}
Let $(r,s,b),(r',s',b')\in\mathfrak{G}$,
$\om'=s'\mathfrak{G}_\nu\in\mathfrak{G}/\mathfrak{G}_\nu$ and
$s'=s'_{\om'}\eta_{\om' s'}$. Then
\begin{eqnarray}
\nonumber (r,s,b)^{-1}(r',s',b')&=&(\f1r,s^{-1},-\f1r s^{-1}bs)(r',s',b')\\
\nonumber&=&(\f{r'}r,s^{-1}s',\f1r s^{-1}(b'-b)s).
\end{eqnarray}
By (\ref{left}) and (\ref{firsttrans'}), we have
\begin{equation}\label{pif}
  f(r',s',b') = e^{-i2\pi \f1{r' }{{s'}^{-1}b's'}\cdot{\nu}}\eta_{\om' s'}^{-1}\a(r',{s'}_{\om'},0)(1-ie_1e_2),
\end{equation}
and
\begin{eqnarray}\label{pi1}
\nonumber  \s_{(r,s,b)}f(r',s',b')&=&f((r,s,b)^{-1}(r',s',b')\\
   &=&e^{-i2\pi \f1{r'} {s'}^{-1}(b'-b)s'\cdot{\nu}}\eta_{\om'' s''}^{-1}\a(\f{r'}r,{s''}_{\om''},0)(1-ie_1e_2),
\end{eqnarray}
where $s''=s^{-1}s'$ and
$\om''\sim s''\mathfrak{G}_\nu$.

Then (\ref{pi1}) is an explicit representation of $\mathfrak{G}$ over
$\mathcal{F}$.

\section{equivalent representations}
In the above section we obtain $\s$, the induced representation of $\mathfrak{G}$. But we can not use it directly. We need alternative versions of $\s$ in order to deal with the symmetry of the Hilbert transformation.

\subsection{}
 Define a linear mapping from $\mathcal{F}$ to $L^2(\RR^3)$ by
\begin{equation}\label{map-FL}
  (\Phi f)(\ul{x})=\a(|\ul{x}|,s_{\f{\ul{x}}{|\ul{x}|}},0),
\end{equation}
where $f$ is represented by (\ref{pif}). Then we obtain  a representation $\Phi^*\s$ of  $\mathfrak{G}$ over $L^2(\RR^3)$.
 Although not a faithful representation,  $\Phi^*\s$ is irreducible. Otherwise, if it was the direct sum of nontrivial representations $\t_1,\t_2$ over $X_1$ and $X_2$ respectively, where $L^2(\RR^3)=X_1\oplus X_2$, then $\s$ should be the direct sum of  $\Phi_*\t_1$ and  $\Phi_*\t_2$.

In the case of dimension-$1$ we know that it is essentially the fact that the Hilbert transformation admits the symmetry of the \textbf{ax+b} group. But by Theorem \ref{schur}  there does not exist a non trivial linear map from $L^2(\RR^3)$ to itself which admits the symmetry of $\mathfrak{G}$.

\subsection{}
To obtain more useful representations we next derive more concise forms of (\ref{pi1}). At first we work on the quaternionic field  $\HH$ instead of $\mathrm{Cl}_3$.

 Let $e_1,e_2,e_3$ be the generators of $\HH$ such that $e_3=e_1e_2$. We define
\begin{equation}\label{}
 \nonumber  S'^+_2:=\CC(1-ie_3)=\{\lm(1-ie_3):\lm\in\CC\},
  \end{equation}
  \begin{equation}\label{}
 \nonumber  S'^-_2:=\CC(1-ie_3)=\{\lm(e_1+ie_2):\lm\in\CC\}
  \end{equation}
  and
  \begin{equation}\label{}
 \nonumber  S'_2:=S'^+_2\oplus S'^-_2.
  \end{equation}
 It is obvious that $S'_2$ is the left ideal of $\HH$ and $S_2'^\pm$ are  linearly isomorphic to  $S^\pm$ respectively.  $\rho^+$ still denotes the irreducible representation of $\textrm{Spin}(2)$ over $S_2'^+$. So we still use $S_2^+$ instead of $S_2'^+$ at the moment. Then the element $f$ of $\mathcal{F}_0$ may be rewritten as
\begin{equation}\label{}
\nonumber  f(x)=\a(x)(1-ie_3), \ x\in\mathfrak{G},
\end{equation}
where $\a(x)$ is a scalar function. Then $f$ satisfies the relation
\begin{equation}\label{}
 \nonumber f((r,s,b)(1,\eta,c)) =e^{-i2\pi c\cdot {\nu}}\eta^{-1}\a(r,s,b)(1-ie_3).
\end{equation}


We consider $\a(r',s'_{\om'},0)$, fixing $r'$
as a function of $\om'$. It  is still  denoted by
$\a(r',\om,0)$. Recall that $\nu=e_3$ and notice that $s'\eta_{\om' s'}^{-1}=s'_{\om'}$ and ${s'}^{-1}b's'\cdot{\nu}=b'\cdot s'\nu{s'}^{-1}=b'\cdot\om'$. The relations (\ref{pif}) and (\ref{pi1}) may be together  rewritten as
 \begin{eqnarray}\label{pi2}
   f(r',s',b') &=& e^{-i2\pi \f1{r' }{b'}\cdot{\om}}{s'}^{-1}s'_{\om'}\a(r',\om',0)(1-ie_3) \\
  \nonumber \s_{(r,s,b)}f(r',s',b') &=&e^{-i2\pi \f1{r'} (b'-b)\cdot \om}(s^{-1}s')^{-1}s''_{\om''}\a(\f{r'}r, \om'',0)(1-ie_3).
\end{eqnarray}
Denote by $s_\om=\gamma(\om)$  the function of $\om$. Noticing that $\om''=s'^{-1}\om's'$,
now (\ref{pi2}) can be rewritten as
\begin{eqnarray}\label{pi3}
  \nonumber f(y) &=& e^{-i2\pi \f1{r' }{b'}\cdot{\om'}}{s'}^{-1}\gamma(\om')\a(r',\om',0)(1-ie_3) \\
 \s_xf(y) &=&e^{-i2\pi \f1{r'} (b'-b)\cdot \om'}{s'}^{-1}s\gamma(s^{-1}\om' s)\a(\f{r'}r, s^{-1}\om' s,0)(1-ie_3).
\end{eqnarray}
For $y=(r',s',b')\in\mathfrak{G}$ let $\ul{x}=r'\om'$ and
  $$f_0(\ul{x})=\gamma(\om)\a(r',\om',0)(1-ie_3).$$
Recall that the invariant measure of $\mathfrak{G}/\mathfrak{G}_\nu$ is $d\mu(y\mathfrak{G}_\nu)=\f1{r'}dr'd\s(\om).$ The correspondence $f\rt f_0$ is a unitary map from $\mathcal{F}$
to the subspace $\tilde{H}^-_1$ of $L^2(\RR^3, \RR_{0,3},\f1{\|x\|^3}d\ul{x})$
whose inverse is given by $f(y)=e^{-i2\pi \f1{r'
}{b'}\cdot{\om'}}{s'}^{-1}f_0(\ul{x})s'$. Notice that
\begin{equation}\label{pi7}
\tilde{H}^-_1=\{\gamma(\om')\a(r',\om',0)(1-ie_3): r'\om'\in \RR^3,
\a(r',\om,0)\in L^2(\RR^3)\}.
\end{equation}

If we conjugate $\s$ by $f\rt
f_0$, by (\ref{pi3}) we obtain an equivalent representation $\pi'$
on $\tilde{H}^-_1$, namely
\begin{equation}\label{pi4}
    \pi'_{(r,s,b)}f_0(\ul{x})=e^{i2\pi \f1{\|\ul{x}\|^2} b\cdot \ul{x}}sf_0(\f1{r}s^{-1}\ul{x}s).
\end{equation}

We now make a change of variable
\begin{equation}
 \nonumber \ul{y}=-\ul{x}^{-1}=\f {\ul{x}}{\|x\|^2},\ g(y)=\|y\|^{-\f32}f_0(\f {\ul{y}}{\|y\|^2}).
\end{equation}
Then the map $f_0\rt g$ isometrically transforms $\tilde{H}^-_1$ onto itself. By (\ref{pi4}),
this map translates $\pi'$ into another representation
\begin{equation}\label{pi5}
    \tilde{\pi}^-_{(r,s,b)}g(\ul{y})=e^{i2\pi  \ul{b}\cdot \ul{y}}r^{\f32}sg({r}s\ul{y}s^{-1}), \ g\in \tilde{H}^-_1.
\end{equation}

\vskip3mm
 On the other hand,  the representation $\rho^-$ of
$\mathrm{Spin}(2)$ on $S^-_2$ induces a representation  of
$\RR^3\times(\RR^+\oplus\mathrm{Spin}(3))$ by the same method  as used above.
  So (\ref{firsttrans}) should be translated into the form
  \begin{equation}\label{firsttrans-}
 f((r,s,b)(1,\eta,c)) =e^{-i2\pi c\cdot {\nu}}\eta^{-1}\a(r,s,b)e_1(1-ie_3),
\end{equation}
and (\ref{pi3}) takes the form
\begin{eqnarray}\label{pi3-}
  \nonumber f(y) &=& e^{-i2\pi \f1{r' }{b'}\cdot{\om'}}{s'}^{-1}\gamma(\om)\a(r',\om',0)(e_1+ie_2)s' \\
 \s_xf(y) &=&e^{-i2\pi \f1{r'} (b'-b)\cdot \om'}s'^{-1}s\gamma(s^{-1}\om' s)\a(\f{r'}r, s^{-1}\om' s,0)(e_1+ie_2),
\end{eqnarray}
which is equivalent with
\begin{equation}\label{pi6}
    \tilde{\pi}^+_{(r,s,b)}g(\ul{y})=e^{i2\pi  \ul{b}\cdot \ul{y}}r^{\f32}sg({r}s\ul{y}s^{-1}), \ g\in
    \tilde{H}^+_1,
\end{equation}
where
\begin{equation}\label{pi8}
\tilde{H}^+_1:=\{\gamma(\om')\a(r',\om',0)(e_1+ie_2): r'\om'\in \RR^3,
\a(r',\om',0)\in L^2(\RR^3)\}.
\end{equation}
\vspace{2mm}

\noindent We summarize (\ref{pi7} ), (\ref{pi5} ), (\ref{pi6} ) into the following theorem.
\bT\label{m0}
The group $r\mathrm{Spin}(3)+\RR^3$ in the quaternionic setting has exactly two irreducible representations read as \\
\indent {\rm(1)}
\begin{equation}\label{pi5'}
    \tilde{\pi}^-_{(r,s,b)}g(\ul{y})=e^{i2\pi  \ul{b}\cdot \ul{y}}r^{\f32}sg({r}s\ul{y}s^{-1}), \ g\in \tilde{H}^-_1,\quad (r,s,b)\in r\mathrm{Spin}(3)+\RR^3,
\end{equation}
where
\begin{equation}\label{pi7'}
\tilde{H}^-_1=\{g(\ul{y})=\gamma(\f{\ul{y}}{|\ul{y}|})\a(|\ul{y}|,\f{\ul{y}}{|\ul{y}|},0)(1-ie_3): \ul{y}\in \RR^3,
\a(|\ul{y}|,\f{\ul{y}}{|\ul{y}|},0)\in L^2(\RR^3)\};
\end{equation}
and
{\rm(2)}
\begin{equation}\label{pi6'}
    \tilde{\pi}^+_{(r,s,b)}g(\ul{y})=e^{i2\pi  \ul{b}\cdot \ul{y}}r^{\f32}sg({r}s\ul{y}s^{-1}), \ g\in
    \tilde{H}^+_1,\quad (r,s,b)\in r\mathrm{Spin}(3)+\RR^3,
\end{equation}
where
\begin{equation}\label{pi8'}
\tilde{H}^+_1=\{g(\ul{y})=\gamma(\f{\ul{y}}{|\ul{y}|})\beta(|\ul{y}|,\f{\ul{y}}{|\ul{y}|},0)(e_1+ie_2): \ul{y}\in \RR^3,
\beta(|\ul{y}|,\f{\ul{y}}{|\ul{y}|},0)\in L^2(\RR^3)\}.
\end{equation}
\eT
Notice that $S'_2e_1=S'_2e_2=\{se_1,s\in S'_2\}$ is also a left ideal of the quaternionic field $\HH$. We can, in the above argument, replace $S'^\pm_2$ with, respectively,
$$S'^+_2e_1=\{\lambda (e_1-e_2),\lambda\in \CC\}$$ and $$S'^-_2e_1=\{\lambda(1+ie_3),\lambda\in \CC\}.$$
That is, the function space ${\mathcal F}_0$ in (\ref{F0}) is changed into
 \begin{equation}
 \begin{array}{rl}
 \{f\in \mathrm{C}(\mathfrak{G}, S'^+_2e_1):\quad&\text{the support of}\ f\ \text{is compact and}\\[1mm]
 &f(x\z)=\rho(\z^{-1})f(x),\ \forall x\in \mathfrak{G},\ \z\in \mathfrak{G}_\nu\},\end{array}
 \end{equation}
 Repeating the argument to prove Theorem \ref{m0},
 we obtain that $\tilde{\pi}$ can be given equivalently as follows.
\bT\label{m00}
The group $r\mathrm{Spin}(3)+\RR^3$ has exactly two irreducible representations given by \\
\indent {\rm(1)}
\begin{equation}\label{pi5''}
    \tilde{\pi}^-_{(r,s,b)}g(\ul{y})=e^{i2\pi  \ul{b}\cdot \ul{y}}r^{\f32}sg({r}s\ul{y}s^{-1}), \ g\in \tilde{H}^-_2,\quad (r,s,b)\in r\mathrm{Spin}(3)+\RR^3,
\end{equation}
where
\begin{equation}\label{pi7''}
\tilde{H}^-_2:=\{g(\ul{y})=\gamma(\f{\ul{y}}{|\ul{y}|})\a(|\ul{y}|,\f{\ul{y}}{|\ul{y}|},0)(e_1-ie_2): \ul{y}\in \RR^3,
\a(|\ul{y}|,\f{\ul{y}}{|\ul{y}|},0)\in L^2(\RR^3)\};
\end{equation}
and\\
{\rm(2)}
\begin{equation}\label{pi6''}
    \tilde{\pi}^+_{(r,s,b)}g(\ul{y})=e^{i2\pi  \ul{b}\cdot \ul{y}}r^{\f32}sg({r}s\ul{y}s^{-1}), \ g\in
    \tilde{H}^+_2,\quad (r,s,b)\in r\mathrm{Spin}(3)+\RR^3,
\end{equation}
where
\begin{equation}\label{pi8''}
\tilde{H}^+_2:=\{g(\ul{y})=\gamma(\f{\ul{y}}{|\ul{y}|})\beta(|\ul{y}|,\f{\ul{y}}{|\ul{y}|},0)(1+ie_3): \ul{y}\in \RR^3,
\beta(|\ul{y}|,\f{\ul{y}}{|\ul{y}|},0)\in L^2(\RR^3)\}.
\end{equation}
\eT
\vskip2mm
\noindent It can be checked that the representation $\tilde{\pi}$ over $\tilde{H}^\pm_1$ are equivalent with those over $\tilde{H}^\pm_2$ respectively because the left multiplications of $\mathrm{Spin}(2)$ over $S^\pm_2$ and $S'^\pm_2e_1$ are identical with each other. In fact, denote by $M_{e_1}$ the right multiplying  by $e_1$ over $S'^+_2$. Then $M_{e_1}$ is a linear isometric isomorphism between the two left ideas $S^\pm_2$ and $S'^\pm_2e_1$. Thus $\rho^\pm M_{e_1}=M_{e_1}\rho^\pm$.

\subsection{} Now we derive two equivalent representations of $\mathfrak{G}$ over the space $L^2(\RR^3,\mathrm{CL}_3)$.
There are some inconveniences to deal with the representations because $S_2$ is not the left ideal of $\mathrm{CL}_3$.
 We revise $S^\pm$ as follows. Let
 \begin{equation}\label{w2}
   W^+_2:=\{w=(1-e_3e_1e_2)\eta,\eta\in S^+_2\}
 \end{equation}
 and
 \begin{equation}\label{w2-}
   W^-_2:=\{w=(1-e_3e_1e_2)\eta,\eta\in S^-_2\}.
 \end{equation}
Define a representation of $\mathrm{Spin}(2)$ by
\begin{equation}\label{w-rep}
  (\rho^+_W s)(w)=sw,\quad \text{for} \ s\in \mathrm{Spin}(2),\ w\in  W^+_2,
\end{equation}
and
\begin{equation}\label{w-rep-}
  (\rho^-_W s)(w)=sw,\quad \text{for} \ s\in \mathrm{Spin}(2),\ w\in  W^-_2.
\end{equation}

Denote by $M_1$ the left multiplication by $1-e_3e_1e_2$ on $S_2$. It is easy to check that $M_1$ is an isometric isomorphisms from, respectively, $S^\pm_2$ to $W^\pm_2$. Recall that elements of $\mathrm{Spin}(3)$ are bi-vectors. Then $s(1-e_3e_1e_2)=(1-e_3e_1e_2)s$ for every $s\in\mathrm{Spin}(3)$. So, $(\rho^+_W s)M_1=M_1(\rho^+_W s)$. We thus  obtain the following proposition.
\bP $\rho^\pm_W$ are respectively irreducible representation of $\mathrm{Spin}(2)$ over $W^\pm_2$, and they are  equivalent with $\rho^\pm$.\eP

By Lemma \ref{gv} the irreducible representations of $\mathfrak{G}_\nu$ over $W^+_2$ are of the form
\begin{equation}\label{}
    \rho_W(s,b)x=e^{i2\pi b\cdot {\nu}}sx,\ \text{for} \ (s,b)\in\mathfrak{G}_\nu, x\in W^+_2.
\end{equation}

As in the previous section,  we consider
 the function space
\begin{equation}\label{F1}
   \begin{array}{rl}
    {\mathcal F}_1
    =\{f\in \mathrm{C}(\mathfrak{G},W^+_2):\quad &\text{the support of }\ f\ \text{is compact  and}\\[1mm]
     &f(x\z)=\rho_W(\z^{-1})f(x),\ \forall x\in \mathfrak{G},\ \z\in \mathfrak{G}_\nu\},
     \end{array}
\end{equation}
where $x=(r,s,b)\in r\mathrm{Spin}(3))+\RR^3$,
$\z=(1,\eta,c)\in\mathrm{Spin}(2)+\RR^3$ and
$\rho_W(\z^{-1})=\rho^{-1}_W(\z)=e^{-i2\pi \nu\cdot
{c}}(\rho^+_W({\eta}))^{-1}$. Let $F'$ be the Hilbert space completion of $F_1$.
  Thus we  get a unitary
representation $\s'$ of $\mathfrak{G}$ on ${\mathcal F}'$ by the left
translation $L'_{(r,s,b)},\ (r,s,b)\in G$,
\begin{equation}\label{}
    L'_{(r,s,b)}(f)(r',s',b')=f({(r,s,b)}^{-1}(r',s',b')),\ \forall {(r',s',b')}\in \mathfrak{G},
\end{equation}
that is, $\s'(x)=L'_x$, being an induced representation
by $\rho_W,$ denoted by
$\mathbf{Ind}^\mathfrak{G}_{\mathfrak{G}_\nu}(\rho_W)$.

Adopting the process  of subsection 4.2 step by step. We obtain the following theorem.
\bT\label{m1}
The group $r\mathrm{Spin}(3)+\RR^3$ has exactly two irreducible representations $\tilde{\pi}'^\pm$ as follows: \\
\indent {\rm(1)}
\begin{equation}\label{pi55}
\nonumber    \tilde{\pi}'^-_{(r,s,b)}g(\ul{y})=e^{i2\pi  \ul{b}\cdot \ul{y}}r^{\f32}sg({r}s\ul{y}s^{-1}), \ g\in \tilde{H}'^-_1,\quad (r,s,b)\in r\mathrm{Spin}(3)+\RR^3,
\end{equation}
where
\begin{equation}\label{pi77}
\nonumber \tilde{H}'^-_1:=\{g(\ul{y})=\gamma(\f{\ul{y}}{|\ul{y}|})\b(\ul{y})(1-e_3e_1e_2)(1-ie_1e_2): \ul{y}\in \RR^3,
\b(\ul{y})\in L^2(\RR^3)\};
\end{equation}
{\rm(2)}
\begin{equation}\label{pi66}
 \nonumber   \tilde{\pi}'^+_{(r,s,b)}g(\ul{y})=e^{i2\pi  \ul{b}\cdot \ul{y}}r^{\f32}sg({r}s\ul{y}s^{-1}), \ g\in
    \tilde{H}'^+_1,\quad (r,s,b)\in r\mathrm{Spin}(3)+\RR^3,
\end{equation}
where
\begin{equation}\label{pi88}
\nonumber \tilde{H}'^+_1:=\{g(\ul{y})=\gamma(\f{\ul{y}}{|\ul{y}|})\beta(\ul{y})(1-e_3e_1e_2)(e_1+ie_2): \ul{y}\in \RR^3,
\beta(\ul{y})\in L^2(\RR^3)\}.
\end{equation}
\eT

Denote by
\begin{eqnarray}
  \nonumber W^+_2e_1 &:=& \{\lambda(1-e_3e_1e_2)(1-ie_1e_2)e_1,\lambda\in\CC\} \\
 \nonumber   &=& \{\lambda(1-e_3e_1e_2)(e_1-ie_2),\lambda\in\CC\}, \\
 \nonumber W^-_2e_1 &:=& \{\lambda(1-e_3e_1e_2)(e_1+ie_2)e_3,\lambda\in\CC\} \\
 \nonumber&=& \{\lambda(1-e_3e_1e_2)(1+ie_1e_2),\lambda\in\CC\}; \\
 \nonumber W^+_2e_3 &:=& \{\lambda(1-e_3e_1e_2)(1-ie_1e_2)e_3,\lambda\in\CC\}, \\
  \nonumber W^-_2e_3 &:=& \{\lambda(1-e_3e_1e_2)(e_1+ie_2)e_3,\lambda\in\CC\}; \\
 \nonumber  W^+_2e_1e_3 &:=& \{\lambda(1-e_3e_1e_2)(1-ie_1e_2)e_1e_3,\lambda\in\CC\} \\
 \nonumber&=& \{\lambda(1-e_3e_1e_2)(e_1-ie_2)e_3,\lambda\in\CC\}, \\
 \nonumber W^-_2e_1e_3 &:=& \{\lambda(1-e_3e_1e_2)(e_1+ie_2)e_1e_3,\lambda\in\CC\}\\
 \nonumber&=& \{\lambda(1-e_3e_1e_2)(1+ie_1e_2)e_3,\lambda\in\CC\}.
\end{eqnarray}

As the same as the above section, we can continue to define the equivalent representations of $\rho^\pm_W$ over $W^\pm_2e_1,W^\pm_2e_3$ and $W^\pm_2e_1e_3$. Then they induce three pairs of representations. We first define six Hilbert spaces:
\begin{eqnarray}
 \nonumber \tilde{H}'^-_2&:=&\{g(\ul{y})=\gamma(\f{\ul{y}}{|\ul{y}|})\b(\ul{y})(1-e_3e_1e_2)(e_1-ie_2): \ul{y}\in \RR^3,
\b(\ul{y})\in L^2(\RR^3)\};\\
 \nonumber \tilde{H}'^+_2&:=&\{g(\ul{y})=\gamma(\f{\ul{y}}{|\ul{y}|})\beta(\ul{y})(1-e_3e_1e_2)(1+ie_1e_2): \ul{y}\in \RR^3,
\beta(\ul{y})\in L^2(\RR^3)\};\\
 \nonumber \tilde{H}'^-_3&:=&\{g(\ul{y})=\gamma(\f{\ul{y}}{|\ul{y}|})\b(\ul{y})(1-e_3e_1e_2)(1-ie_1e_2)e_3: \ul{y}\in \RR^3,
\b(\ul{y})\in L^2(\RR^3)\};\\
 \nonumber \tilde{H}'^+_3&:=&\{g(\ul{y})=\gamma(\f{\ul{y}}{|\ul{y}|})\beta(\ul{y})(1-e_3e_1e_2)(e_1+ie_2)e_3: \ul{y}\in \RR^3,
\beta(\ul{y})\in L^2(\RR^3)\};\\
 \nonumber \tilde{H}'^-_4&:=&\{g(\ul{y})=\gamma(\f{\ul{y}}{|\ul{y}|})\b(\ul{y})(1-e_3e_1e_2)(e_1-ie_2)e_3: \ul{y}\in \RR^3,
\b(\ul{y})\in L^2(\RR^3)\};\\
 \tilde{H}'^+_4&:=&\{g(\ul{y})=\gamma(\f{\ul{y}}{|\ul{y}|})\beta(\ul{y})(1-e_3e_1e_2)(1+ie_1e_2)e_3: \ul{y}\in \RR^3,
\beta(\ul{y})\in L^2(\RR^3)\}.
\end{eqnarray}

\bT\label{m11}
The group $r\mathrm{Spin}(3)+\RR^3$ has the following irreducible representations\\
\indent {\rm(1)}
\begin{equation}\label{pi551}
    \tilde{\pi}'^\pm_{(r,s,b)}g(\ul{y})=e^{i2\pi  \ul{b}\cdot \ul{y}}r^{\f32}sg({r}s\ul{y}s^{-1}), \quad (r,s,b)\in r\mathrm{Spin}(3)+\RR^3,\  g\in \tilde{H}'^\pm_j,\ j=2,3,4.
\end{equation}
Moreover, each of them is equivalent to its corresponding representation over $ H'^\pm_1.$
\eT

\subsection{}Now we present the results in $\RR^2$ with the $\mathrm{Cl}_2$ setting.

Define
 \begin{eqnarray}
 U^+_2 &:=& \{\lambda(e_1+ie_2+1-ie_1e_2),\lambda\in\CC\} \\
 U^-_2&:=& \{\mu(e_1+ie_2-(1-ie_1e_2)),\lambda\in\CC\}.
\end{eqnarray}
Denote by $\rho^\pm$ the left multiplications: $\rho^\pm(s)x=sx$ for $x\in U^\pm_2.$
 It is obvious that they are all equivalent to each other. Those representations can induce the representations of $\mathfrak{G}$ for characterizing the Hilbert transformation.

So, up to equivalence, there exists
only one irreducible representation ${\pi'}^+$ of $G_\nu$ such that
\begin{equation}\label{2repre1}
  ({\pi'}^+(1,s,b))x=e^{i2\pi<b,(0,1)>}sx,\ s\in\mathrm{spin}(1)=\ZZ_2, x\in U^+_2.
\end{equation}

  By the Mackey machine \cite{ka}, we can induce a
representation of $\mathfrak{G}$ by $\rho$ , for which we denote by
${\mathcal C}(\mathfrak{G},U^+_2)$ the space of continuous Weyl
spinor-valued functions from $\mathfrak{G}$ to $U^+_2$. Then we
begin with a functions space as follows:
\begin{equation}\label{2F0}
   \begin{array}{rl}
    {\mathcal F}_0
    =\{f\in {\mathcal C}(\mathfrak{G},U^+_2):\quad &\text{the support of }\ f\ \text{is compact  and}\\[1mm]
     &f(x\z)=\rho(\z^{-1})f(x),\ \forall x\in \mathfrak{G},\ \z\in \mathfrak{G}_\nu\},
     \end{array}
\end{equation}
where $x=(r,s,b)\in r\mathrm{Spin}(2))+\RR^2$,
$\z=(1,\eta,c)\in\mathrm{Spin}(1)+\RR^2$ and
$\rho(\z^{-1})=\rho^{-1}(\z)=e^{-i2\pi \nu\cdot
{c}}(\rho^+({\eta}))^{-1}$. Recalling that $\rho^+(\eta)$ is just the left multiplication of $\eta$, thus the equality in (\ref{2F0}) may be
rewritten as
\begin{equation}\label{2coset}
  f((r,s,b)(1,\eta,c))=e^{-i2\pi {c}\cdot\nu}\eta^{-1}f(r,s,b).
\end{equation}
Notice that $\eta$ is either $1$ or $-1$. So by (\ref{2coset}) $f(r,s,b)=f(r,-s,b)$.

Denoted by $<\cdot,\cdot>_0$ the inner product of the complexification
of $\RR_{0,2}$. For $f,g\in{\mathcal F}_0$,
$<f(x),g(x)>_0$ depends only on the $\mathfrak{G}_\nu$-coset of $x$.
So $<f(x),g(x)>_0$ can be considered as a complex function on
$\mathfrak{G}/\mathfrak{G}_\nu$. Thus we can define an inner product
on ${\mathcal F}_0$ such that
\begin{equation}\label{2product}
    <f,g>=\int_{\mathfrak{G}/\mathfrak{G}_\nu}<f(x),g(x)>_0 d\mu(x \mathfrak{G}_\nu).
\end{equation}

 Let ${\mathcal F}$ be the
Hilbert space completion of ${\mathcal F}_0$. So we  get a unitary
representation $\s$ of $\mathfrak{G}$ on ${\mathcal F}$ by the left
translation $L_x,\ x\in G$,
\begin{equation}\label{2left}
    L_x(f)(y)=f(x^{-1}y),\ \forall y\in \mathfrak{G},
\end{equation}
that is, $\s(x)=L_x$. Then $\s$ is called the induced representation
by $\rho$ and denoted by
$\mathbf{Ind}^\mathfrak{G}_{\mathfrak{G}_\nu}(\rho)$.
\vspace{3mm}

\noindent We need to parameterize the elements of $\mathfrak{G}/\mathfrak{G}_\nu$
and those of ${\mathcal F}$. Since
\begin{equation}\label{}
\nonumber\mathfrak{G}/\mathfrak{G}_\nu=(\RR^2\times(\RR^+\oplus\mathrm{Spin}(2)))/(\mathrm{\RR^2\ltimes Spin}(1)),
\end{equation}
we obtain
\begin{equation}\label{2paramet}
\mathfrak{G}/\mathfrak{G}_\nu\cong \RR^+\times
(\mathrm{Spin}(2)/\mathrm{Spin}(1))\cong \RR^+\times S^1,
\end{equation}
where $\cong$ means the isometric homeomorphism.

By (\ref{2F0}) and (\ref{2coset}) we can define a function
$\tilde{f}(\tilde{x})$ over $\mathfrak{G}/\mathfrak{G}_\nu$ by
$\tilde{f}(\tilde{x})=f(x)$ for $\tilde{x}=x\mathfrak{G}_\nu\in
\mathfrak{G}/\mathfrak{G}_\nu$. This means that the functions $f$ are essentially those
 over $\mathfrak{G}/\mathfrak{G}_\nu$.

The element $f$ of $\mathcal{F}_0$ may be rewritten as
\begin{equation}\label{}
\nonumber  f(x)=\a(x)(e_1+ie_2+1+ie_1e_2), x\in\mathfrak{G},
\end{equation}
where $\a(x)$ is a scalar function. Then $f$ satisfies
\begin{equation}\label{2firsttrans}
 f((r,s,b)(1,\eta,c)) =e^{-i2\pi c\cdot {\nu}}\eta^{-1}\a(r,s,b)(e_1+ie_2+1+ie_1e_2).
\end{equation}

Recall that  there exists an isometrical homeomorphism
between   $\mathfrak{G}/\mathfrak{G}_\nu$ and the unit circle  $S^1$ by the correspondence $s\rightarrow \om =s\nu s^{-1}$.
For $\om\in S^1$, choose $s_\om\in \mathfrak{G}$ such that  $\tilde{s}_\om=s_\om
\mathfrak{G}_\nu$ is its corresponding element in
$\mathfrak{G}/\mathfrak{G}_\nu$, i.e., ${\om}=s_\om e_2s^{-1}_\om$. Let $s_\om=\cos(\phi)+\sin(\phi)e_1e_2$. Then
$$
\om(\cos(\phi)+\sin(\phi)e_1e_2)=(\cos(\phi)+\sin(\phi)e_1e_2)e_2.
$$
We always choose $\phi\in[0,\pi)$.

Recall that $(r,s,b)(1,\om,c)=(r,s\om,rscs^{-1}+b)$. If we take
$b=0,s=s_\om,\eta=\eta_{\om s}$ and
write $\f1r s^{-1}bs$ in place of $c$ in (\ref{2firsttrans}), we have
\begin{eqnarray}\label{2firsttrans'}
 f(r,s,b) &=&f((r,s_\om,0)(1,\eta_{\om s},\f1r s^{-1}bs)\\
 \nonumber&=& e^{-i2\pi \f1r {s^{-1}bs}\cdot{\nu}}\eta_{\om s}^{-1}\a(r,s_\om,0)(e_1+ie_2+1+ie_1e_2).
\end{eqnarray}
Let $x=(r,s,b),y=(r',s',b')\in\mathfrak{G}$,
$\om'=s'\mathfrak{G}_\nu\in\mathfrak{G}/\mathfrak{G}_\nu$ and
$s'=s'_{\om'}\eta_{\om' s'}$. Then
by (\ref{2left}) and (\ref{2firsttrans'}), we have
\begin{eqnarray}\label{2pi1}
 \nonumber f(y) &=& e^{-i2\pi \f1{r' }{{s'}^{-1}b's'}\cdot{\nu}}\eta_{\om' s'}^{-1}\a(r',{s'}_{\om'},0)(e_1+ie_2+1+ie_1e_2), \\
  \s_xf(y)=f(x^{-1}y)&=&e^{-i2\pi \f1{r'} {s'}^{-1}(b'-b)s'\cdot{\nu}}\eta_{\om'' s''}^{-1}\a(\f{r'}r,{s''}_{\om''},0)(e_1+ie_2+1+ie_1e_2),
\end{eqnarray}
where $s''=s^{-1}s'$ and
$\tilde{\om''}=s''\mathfrak{G}_\nu$.

Then (\ref{2pi1}) is the explicit representation of $\mathfrak{G}$ over
$\mathcal{F}$.
\vspace{2mm}

Denote by
\begin{equation}\label{2pi7}
\tilde{\tilde{H}}^-_1:=\{s_{\om'}\a(r',\om',0)(e_1+ie_2+1+ie_1e_2): r'\om'\in \RR^3,
\a(r',\om,0)\in L^2(\RR^3)\}
\end{equation}
and
\begin{equation}\label{2pi8}
\tilde{\tilde{H}}^+_1:=\{s_{\om'}\a(r',\om',0)(e_1+ie_2-(1-ie_1e_2)): r'\om'\in \RR^2,
\a(r',\om',0)\in L^2(\RR^2)\}.
\end{equation}

\noindent We get the  the following theorem.
\bT\label{2m0}
The group $r\mathrm{Spin}(2)+\RR^2$ has two  equivalently irreducible representations as follows: \\
\indent {\rm(1)}
\begin{equation}\label{2pi5'}
    \tilde{\tilde{\pi}}^-_{(r,s,b)}g(\ul{y})=e^{i2\pi  \ul{b}\cdot \ul{y}}r^{\f32}sg({r}s\ul{y}s^{-1}), \ g\in \tilde{\tilde{H}}^-_1,\quad (r,s,b)\in r\mathrm{Spin}(2)+\RR^2.
\end{equation}
{\rm(2)}
\begin{equation}\label{2pi6'}
    \tilde{\tilde{\pi}}^+_{(r,s,b)}g(\ul{y})=e^{i2\pi  \ul{b}\cdot \ul{y}}r^{\f32}sg({r}s\ul{y}s^{-1}), \ g\in
    \tilde{\tilde{H}}^+_1,\quad (r,s,b)\in r\mathrm{Spin}(2)+\RR^2.
\end{equation}
\eT

Substitute $U^\pm_2 e_1=\{we_1:\ w\in U^\pm\}$ for $U^\pm_2$ in the argument proving Theorem \ref{2m0}. Setting
\begin{equation}\label{2pi7''}
\tilde{\tilde{H}}^-_2:=\{g(\ul{y})=s_{\f{\ul{y}}{|\ul{y}|}}\a(|\ul{y}|,\f{\ul{y}}{|\ul{y}|},0)(e_1+ie_2+1-ie_1e_2)e_1: \ul{y}\in \RR^2,
\a(|\ul{y}|,\f{\ul{y}}{|\ul{y}|},0)\in L^2(\RR^2)\};
\end{equation}
and
\begin{equation}\label{2pi8''}
\tilde{\tilde{H}}^+_2:=\{g(\ul{y})=s_{\f{\ul{y}}{|\ul{y}|}}\beta(|\ul{y}|,\f{\ul{y}}{|\ul{y}|},0)(e_1+ie_2-(1-ie_1e_2))e_1: \ul{y}\in \RR^2,
\beta(|\ul{y}|,\f{\ul{y}}{|\ul{y}|},0)\in L^2(\RR^2)\},
\end{equation}
we obtain
\bT
The group $r\mathrm{Spin}(2)+\RR^2$ has the following two equivalently irreducible representations:
\begin{equation}\label{2pi5''}
    \tilde{\tilde{\pi}}^\pm_{(r,s,b)}g(\ul{y})=e^{i2\pi  \ul{b}\cdot \ul{y}}r^{\f32}sg({r}s\ul{y}s^{-1}), \ g\in \tilde{\tilde{H}}^\pm_2,\quad (r,s,b)\in r\mathrm{Spin}(2)+\RR^2.
\end{equation}
\eT

\section{operators with symmetry $\mathfrak{G}$: case of Quaternion}
In section 2  we obtained that both the Hilbert and Riesz transformations
commute with the group $\mathfrak{G}.$ On the
contrary we want to know what happens to the
operators  that possess such symmetry of $\mathfrak{G}$. In the present section we treat this problem by using the quaternions and leave the case of Clifford to the next section.
\subsection{Decomposition of $L^2(\RR^3,S_2)$}

We first need to point out that $L^2(\RR^3)$ is not the right space
to discuss the relation between the Hilbert transformation and
$\mathfrak{G}$. In fact, in (\ref{spinrepn}), we denoted by
$\pi$ the natural representation of $r\mathrm{Spin}(n)+\RR^n$ over
the Clifford-valued square integral functions space
$L^2(\RR^3,\mathrm{Cl}_3)$. For $(r,s,b)\in \mathfrak{G}$,  $\pi
_{rsb}$ is not the map from  $L^2(\RR^3)$ to itself.
$L^2(\RR^3,\mathrm{Cl}_3)$ is also not the right setting. Even  $L^2(\RR^3,\HH)$ is not so since it
is obvious that $L^2(\RR^3,\HH$) can not be the sum of
$H^+_1$ and $H^-_1$. In what follows, we consider the Hilbert transformation
and the naturel representation of $\mathfrak{G}$ in $L^2(\RR^3,S_2)$
where $S_2$ is a minimal left ideal of $\mathrm{Cl}_3$. Of course it is also the minimal left ideal of $\HH$ when we consider $S'_2$ as $S_2$ by taking $e_3=e_1e_2$. In this section we do not make distinction between $S'_2$ and $S_2$.
\vspace{2mm}

\noindent It is tedious but trivial to prove the following proposition.
\bP The Hilbert transformation $H$ and $\pi_{rsb},(r,s,b)\in\mathfrak{G}$ are unitary operators on $L^2(\RR^3,S_2)$. Moreover,
\begin{equation}\label{pp}
  \nonumber (H\pi_{rsb})f=(\pi_{rsb}H)f,\quad f\in L^2(\RR^3)
\end{equation}
if and only if
\begin{equation}\label{}
  \nonumber (H\pi_{rsb})f=(\pi_{rsb}H)f,\quad f\in L^2(\RR^3,S_2).
\end{equation}
\eP

\brs  \rm It is similar to the case of one dimension.  Neither the upper half plane Hardy space, nor the lower half plane Hardy space contain the real-valued function other than $0$. So the Hilbert transformation does not have any proper invariant subspace in $L^2(\RR,\RR),$ and thus the complex-valued
square integrable functions space
$L^2(\RR)=L^2(\RR,\CC)$ is the right setting.
\ers
The following lemma gives the reason of the choice of $L^2(\RR^3,S_2)$.
\bL \label{l2s2}
Hilbert space  $L^2(\RR^3,S_2)$ is the direct sum of $\tilde{H}^+_1$ and $\tilde{H}^-_1$.
\eL

 We first prove a useful lemma.
 \bL \label{de} Let $\chi_\pm$ be defined by (\ref{projection1}).
  \begin{enumerate}
  \item  $f\in \tilde{H}^+_1$ if and only if $f\in L^2(\RR^3,S_2)$ and $f=\chi_+f$;
   \item $f\in \tilde{H}^-_1$ if and only if $f\in L^2(\RR^3,S_2)$ and $f=\chi_-f$.
   \end{enumerate}
\eL

\bpp
We only prove the first statement since the second is similar.
Let $f\in \tilde{H}^+_1$. By (\ref{pi7}) there exists $\beta(|\ul{x}|,\f{\ul{x}}{|\ul{x}|},0)\in L^2(\RR^3)$ such that $f(\ul{x})=\beta(|\ul{x}|,\f{\ul{x}}{|\ul{x}|},0)s_{\f{\ul{x}}{|\ul{x}|}}(e_1+ie_2)$. By the choice of $s_{\f{\ul{x}}{|\ul{x}|}}$ we get $\f{\ul{x}}{|\ul{x}|}=s_{\f{\ul{x}}{|\ul{x}|}}e_3s_{\f{\ul{x}}{|\ul{x}|}}^{-1}$, and
$$
\f{\ul{x}}{|\ul{x}|}f(\ul{x})=s_{\f{\ul{x}}{|\ul{x}|}}e_3s_{\f{\ul{x}}{|\ul{x}|}}^{-1}\beta(|\ul{x}|,\f{\ul{x}}{|\ul{x}|},0)s_{\f{\ul{x}}{|\ul{x}|}}(e_1+ie_2).
$$

Then by $e_3(e_1+ie_2)=e_2-ie_1=-i(e_1+ie_2)$ we have $\f{\ul{x}}{|\ul{x}|}f(\ul{x})=-if(\ul{x})$, which gives
$$
f(\ul{x})=i\f{\ul{x}}{|\ul{x}|}f(\ul{x})=\f12(1+i\f{\ul{x}}{|\ul{x}|})f(\ul{x})=\chi_+(\ul{x}))f(\ul{x}).
$$
On the contrary,  assume that $f(\ul{x})\in L^2(\RR^3,S_2)$ and $f(\ul{x})=\chi_+(\ul{x})f(\ul{x}).$ Then we have $\f{\ul{x}}{|\ul{x}|}f(\ul{x})=-if(\ul{x})$. Since the left multiplication by
$s_{\f{\ul{x}}{|\ul{x}|}}$ is an isometric isomorphism on   $f\in L^2(\RR^3,S_2)$ and $L^2(\RR^3,S_2)=L^2(\RR^3,S^+_2)\bigoplus L^2(\RR^3,S^-_2)$, there exist $\a,\beta\in L^2(\RR^3)$ such that
\begin{equation}\label{decoms2}
  f(\ul{x})=s_{\f{\ul{x}}{|\ul{x}|}}(\a(|\ul{x}|,\f{\ul{x}}{|\ul{x}|},0)(1-ie_3)+\beta(|\ul{x}|,\f{\ul{x}}{|\ul{x}|},0)(e_1+ie_2)).
\end{equation}
Left multiplying ${\f{\ul{x}}{|\ul{x}|}}=s_{\f{\ul{x}}{|\ul{x}|}}e_3s_{\f{\ul{x}}{|\ul{x}|}}^{-1}$ to the both sides of (\ref{decoms2}), we have
\begin{equation}\label{decoms2'}
 \f{\ul{x}}{|\ul{x}|}f(\ul{x})=s_{\f{\ul{x}}{|\ul{x}|}}e_3(\a(|\ul{x}|,\f{\ul{x}}{|\ul{x}|},0)(1-ie_3)+\beta(|\ul{x}|,\f{\ul{x}}{|\ul{x}|},0)(e_1+ie_2)).
\end{equation}
Recall that $e_3(1-ie_3)=i(1-ie_3)$.  Substituting  $\f{\ul{x}}{|\ul{x}|}f(\ul{x})=-if(\ul{x})$ in the left part, by (\ref{decoms2}) and (\ref{decoms2'}), we get
 \begin{equation}\label{decoms2''}
  \nonumber -i s_{\f{\ul{x}}{|\ul{x}|}}\a(|\ul{x}|,\f{\ul{x}}{|\ul{x}|},0)(1-ie_3)= is_{\f{\ul{x}}{|\ul{x}|}}\a(|\ul{x}|,\f{\ul{x}}{|\ul{x}|},0)(e_1+ie_2),
 \end{equation}
 which implies that $\a(|\ul{x}|,\f{\ul{x}}{|\ul{x}|},0)=0$. Then by (\ref{decoms2}) we obtain $f\in \tilde{H}^+_1$.
\epp
\brs
We must point that the identity $e_3(e_1+ie_2)=e_2-ie_1=-i(e_1+ie_2)$ is essential for   the proof of Lemma \ref{de}. In the following discussion of both cases of  $\mathrm{CL}_2$ and $\mathrm{CL}_3$, we can not find such identity in the Weyl spinor space.
\ers
\vspace{3mm}
\noindent \it Proof of Lemma \ref{l2s2}. \rm We only need to prove that $L^2(\RR^3,S_2)\subseteq \tilde{H}^+_1\bigoplus \tilde{H}^-_1$. Suppose that $f\in L^2(\RR^3,S_2)$. Denote $f^+=\chi_+f,f^-=\chi_-f$. It is easy to verify that
\begin{equation}\label{decoms2'''}
f=f^++f^-,\quad \chi_+f^+=f^+, \quad \chi_-f^-=f^-.
 \end{equation}
Now by (\ref{decoms2'''}) and Lemma \ref{de} we have $f^+\in \tilde{H}^+_1$ and $f^-\in \tilde{H}^-_1$.
\vskip2mm
By the same argument we have the following two lemmas.

 \bL \label{de3} Let $\chi_\pm$ be defined by (\ref{projection1}).
\begin{enumerate}
  \item $f\in \tilde{H}^+_2$ if and only if $f\in L^2(\RR^3,S_2e_1)$ and $f=\chi_+f$;
  \item $f\in \tilde{H}^-_2$ if and only if $f\in L^2(\RR^3,S_2e_1)$ and $f=\chi_-g$.
\end{enumerate}

\eL

\bL \label{l2s23}
The Hilbert space  $L^2(\RR^3,S_2e_1)$ is the direct sum of $\tilde{H}^+_2$ and $\tilde{H}^-_2$.
\eL
Now we obtain
\bT
\begin{equation}\label{decH}
  L^2(\RR^3,\HH)=L^2(\RR^3,S_2)\oplus L^2(\RR^3,S_2e_1)=\tilde{H}^+_1\oplus \tilde{H}^-_1\oplus \tilde{H}^+_2\oplus \tilde{H}^-_2.
\end{equation}
\eT
\subsection{Decomposition of $\pi$ } We first need to derive some results in relation to quaternionic Fourier transform.

For $f\in L^2(\RR^3,\HH)$, the Fourier transformation
 is defined by
\begin{equation}\label{four}
\mathfrak{F}(f)(\ul{\xi}):=\int_{\RR^3 } e^{-i2\pi<\ul{x},\ul{\xi}>}f(\ul{x})d\ul{x}
\end{equation}
and the inverse Fourier transformation is defined by
\begin{equation}\label{four-}
\mathfrak{F}^{-1}(g)(\ul{x}):=\int_{\RR^3 } e^{i2\pi<\ul{x},\ul{\xi}>}g(\ul{\xi})d\ul{\xi}.
\end{equation}

Both the Fourier and inverse Fourier transformations, defined respectively in (\ref{four}) and (\ref{four-}), are isometrically automorphism of $L^2(\RR^3,S_2)$. 
Let $^q\!{H}^\pm_j$ be the Hilbert spaces isometrically isomorphic to $\tilde{H}^\pm_j$, $j=1,2$, as images of the Fourier
transformation mappings, i.e., $^q\!{H}^\pm_j=\mathfrak{F}\tilde{H}^\pm_j$. So $L^2(\RR^3,S_2)$ and $L^2(\RR^3,S_2e_1)$ have other direct sum decompositions
\begin{equation}\label{decompi}
   L^2(\RR^3,S_2)= {}^q\!{H}^+_1\oplus {}^q\!{H}^-_1,\ L^2(\RR^3,S_2e_1)={}^q\!{H}^+_2\oplus {}^q\!{H}^-_2.
\end{equation}
 Moreover, we get  two representations   ${\pi}^\pm$ that are equivalent to
$\tilde{\pi}^\pm$ respectively via the Fourier transformation.
\bT\label{ml}
Let $\pi$ be the natural unitary
representation  of $\mathfrak{G}_3$ over $ L_2(\RR^3,S_2)$.
Then $\pi$ is decomposed into the direct sum of two sub-representations of ${\pi}^\pm$ on $^q\!{H}^\pm_1,$ respectively, i.e.,
\begin{equation}\label{decom}
  \nonumber    \pi={\pi}^+\oplus{\pi}^-.
\end{equation}
\eT
\bpp
For $f\in L_2(\RR^3,S_2)$ and  $(r,s,b)\in\mathfrak{G}_3$, let $\check{f}=\mathfrak{F}^{-1}f$. Then there exist $\check{f}^+_1\in \tilde{H}^+_1$ and $\check{f}^-_1\in \tilde{H}^-_1$ such that $\check{f}=\check{f}^+_1+\check{f}^-_1$. Denote by $f^\pm=\mathfrak{F}\check{f}^\pm_1,$ and $f^\pm$ belong to $^q\!{H}^{\pm},$ all respectively, and $f=f^++f^-$.\\
\indent By Theorem \ref{m0} we have
\begin{equation}\label{pi5'''}
    \tilde{\pi}^+_{(r,s,b)}\check{f}^+_1(\ul{y})=e^{i2\pi  \ul{b}\cdot \ul{y}}r^{\f32}s\check{f}^+_1({r}s\ul{y}s^{-1}), \ul{y}\in\RR^3,
\end{equation}
and
\begin{equation}\label{pi6'''}
    \tilde{\pi}^-_{(r,s,b)}\check{f}^-_1(\ul{y})=e^{i2\pi  \ul{b}\cdot \ul{y}}r^{\f32}s\check{f}^-_1({r}s\ul{y}s^{-1}),\ \ul{y}\in\RR^3.
\end{equation}
By the definition of $\pi$, we have
\begin{equation}\label{spinrepn'}
    \pi_{rsb}(f)(\ul{x})=(\f1r)^{\f 32}sf(\f1r s^{-1}(\ul{x}-b)s).
\end{equation}
Changing variable $\ul{x}=rs\ul{\xi}s^{-1}+b$, by (\ref{four}) and (\ref{spinrepn'}), we have
\begin{eqnarray}\label{four'}
 \nonumber \mathfrak{F}^{-1}(\pi_{rsb}(f^+))(\ul{y}) &=& \int_{\RR^3 } e^{i2\pi<\ul{x},\ul{y}>}(\f1r)^{\f 32}sf^+(\f1r
    s^{-1}(\ul{x}-b)s)d\ul{x} \\
  \nonumber  &=& \int_{\RR^3 } e^{i2\pi<rs\ul{\xi}s^{-1}+b,\ul{y}>}(\f1r)^{\f 32}sf(\ul{\xi})r^nd\ul{\xi} \\
 \nonumber  &=&  e^{i2\pi  <\ul{b}, \ul{y}>}r^{\f32}s\int_{\RR^3 } e^{i2\pi<rs\ul{\xi}s^{-1},\ul{y}>}(\f1r)^{\f 32}f(\ul{\xi})r^nd\ul{\xi} \\
&=&  e^{i2\pi  <\ul{b}, \ul{y}>}r^{\f32}s\check{f}^+_1({r}s\ul{y}s^{-1}).
\end{eqnarray}
Similarly, we have
\begin{equation}\label{four''}
  \mathfrak{F}^{-1}(\pi_{rsb}(f^-))(\ul{y})= e^{i2\pi  <\ul{b}, \ul{y}>}r^{\f32}s\check{f}^-_1({r}s\ul{y}s^{-1}).
\end{equation}
Now by (\ref{pi5'''}), (\ref{pi6'''}),  (\ref{four'}) and (\ref{four''}) we have
\begin{equation}\label{}
 \nonumber \mathfrak{F}^{-1}(\pi_{rsb}(f))(\ul{y})=\mathfrak{F}^{-1}(\pi_{rsb}(f^++f^-))(\ul{y})=(\tilde{\pi}^+\oplus\tilde{\pi}^-)(\check{f})(\ul{y}),
\end{equation}
which means that $\pi^\pm=\mathfrak{F}^{-1}\tilde{\pi}^\pm$ are the irreducible subrepresentations of $\pi$ over $^q\!{H}^\pm_1,$ respectively, and
\begin{equation}\label{}
 \nonumber    \pi={\pi}^+\oplus{\pi}^-.
\end{equation}
\epp

We know that $\pi$ is also the natural unitary
representation  of $\mathfrak{G}_3$ over $ L_2(\RR^3,S_2e_1)$. Similarly we also have the following equivalent  representations of   ${\pi}^\pm$.
\bT\label{m2 }
 $\pi$ is decomposed into the direct sum of the sub-representations of ${\pi}^\pm$ on $^q\!{H}^\pm_2,$ respectively. Moreover, the representations of ${\pi}^\pm$ over $^q\!{H}^\pm_1$ are equivalent to their representations over $^q\!{H}^\pm_2$ respectively.
\eT
In \cite{kq,qm}   a series of results concern eigenvalues, invariant spaces, etc., of the Hilbert transformation and their relations with harmonic analysis.  For instance, they prove
\begin{equation}\label{fou}
\mathfrak{F}^{-1}(Hf)(\ul{\xi})=i\f{\ul{\xi}}{|\ul{\xi}|}\check{f},  f\in L_2(\RR^3,\HH),
\end{equation}
where $\mathfrak{F}^{-1}$ is the inverse Fourier transformation.
Here we give their interpretations in the setting of group representation.

Here we first recall several results in Clifford analysis.  A monogenic function $f$ is a solution for the Dirac equation $Df=0$.
Especially, the so called Cauchy kernel function $\f{\ul{y}-\ul{x}}{|y-x|^{n+1}}$ is  monogenic in $\RR^{0,3}\setminus\{y\}$.

Denote by
$$
^q\!{H}^+:= {}^q\!{H}^+_1\oplus {}^q\!{H}^+_2,\ {}^q\!{H}^-={}^q\!{H}^-_1\oplus{}^q\!{H}^-_2.$$

\bT\label{last} Let $H$ be the Hilbert transformation extended to $L^2(\RR^3,\HH)$ based on its definition over $L^2(\RR^3)$ and the relation (\ref{hilbert3}). Then it has two nontrivial maximal invariant
spaces $ {}^q\!{H}^+$ and $ {}^q\!{H}^-$ spelt as
\begin{equation}\label{}
 \nonumber   L^2(\RR^3,S_2)= {}^q\!{H}^+\oplus {}^q\!{H}^-,
\end{equation}
and
\begin{equation}\label{h-charact}
  Hf=f,\quad \forall f\in  {}^q\!{H}^+;\quad Hg=-g,\quad \forall g\in  {}^q\!{H}^-.
\end{equation}
\eT
\bpp It comes from  Proposition \ref{p1} and \ref{pp} that $H$ commutes with the natural representation of $r\mathrm{Spin}(3)+\RR^3$ over $L^2(\RR^3,S_2)$. So $H$ has the same invariant subspaces as $\pi$, i.e., it is invariant over both $ {}^q\!{H}^+$ and $ {}^q\!{H}^-$. Again $H$ commutes with the irreducible representations  $\pi^\pm$ of   $ \mathfrak{G} $ over $ {}^q\!{H}^\pm_j,j=1,2,$ respectively.
Thus by Theorem \ref{t1} (Schur's Lemma) there exist two complex numbers $\lambda,\lambda',\mu$ and $\mu'$ such that
\begin{equation}\label{new1}
  Hf=\lambda f,\forall f\in  {}^q\!{H}_1^+;\ Hg=\mu g,\forall g\in  {}^q\!{H}_1^-
\end{equation}
and
\begin{equation}\label{new2}
  Hf=\lambda' f,\forall f\in  {}^q\!{H}_2^+;\ Hg=\mu' g,\forall g\in  {}^q\!{H}_2^-.
\end{equation}

Fixed $\ul{y}\in \RR^{0,3},y_0<0$, then $F=\f{\ul{y}-\ul{x}}{|y-x|^{4}}$
is the monogenic function  on the closed upper space (containing the boundary $\RR^3)$ in $\RR^{0,3}$ and $F\in L^2(\RR^3,\HH)$. Let $F=F_1+F_2$, where $F_1\in L^2(\RR^3,S_2)$ and $F_2\in L^2(\RR^3,S_2e_1)$. Denote by $\check{F}$ its Fourier transformation.  Thus by Theorem \ref{Plemelj theorem} we have
\begin{equation}\label{ifi0}
HF=F.
\end{equation}
 Take the inverse Fourier transformation on both sides of (\ref{ifi0}). We have $$\mathfrak{F}^{-1}(HF)=\check{F}.$$
 By the above formula and (\ref{fou}) we have
 $\check{F}=i\f{\ul{\xi}}{|\ul{\xi}|}\check{F}$,  and then $\check{F}=\chi_+(\ul{\xi})\check{F}$, i.e.
 $$
 \check{F}_1+\check{F}_2= \chi_+(\ul{\xi})\check{F}_1+\chi_+(\ul{\xi})\check{F}_2.
  $$
  We obtain
  $$
   \check{F}_1= \chi_+(\ul{\xi})\check{F}_1;\ \check{F}_2=\chi_+(\ul{\xi})\check{F}_2,
   $$
   since $\chi_+(\ul{\xi})\check{F}_1$ is still in $L^2(\RR^3,S_2)$, and $\chi_+(\ul{\xi})\check{F}_2$ is still in $L^2(\RR^3,S_2e_1)$.
  Now by Lemma \ref{de} and Lemma \ref{de3} we get that $\check{F}_1\in \tilde{H}^+_1$ and $\check{F}_2\in \tilde{H}^+_2$.  By the definition of $ {}^q\!{H}^+$, we obtain that $F\in  {}^q\!{H}^+$.
But by (\ref{new1}) and (\ref{new2}) we have
\begin{equation}\label{ifi}
 HF=HF_1+HF_2=\lambda F_1+\lambda'F_2.
 \end{equation}
Comparing (\ref{ifi0}) and (\ref{ifi}) gives $\lambda=\lambda'=1$.

 By the similar procedure we can get that $\mu=\mu'=-1$.
 \epp
 \vspace{2mm}

\noindent  By  (\ref{h-charact}) and the Plemelj formulas (Theorem \ref{Plemelj theorem}), we get the following corollary.

 \bC $H^\pm$ in Theorem \rm\ref{last} \it are
 the upper and lower monogenic Hardy spaces, respectively.\eC

Because  the Lebesgue integral commutes with the singular integral,  $H$ is  self-adjoint over $L^2(\RR^3,\HH)$.  So its point spectra $\lambda,\mu$ must be real numbers.

\vspace{3mm}

\noindent We end this section through proving the inversion of Theorem \ref{last}.
\bT\label{qlast}
Let $T$ be a  linear  operator over $L^2(\RR^3,\HH)$. If  $T$ is  commutes with the natural representation of $\mathfrak{G}$. Then $T$ must be be of the form  $\lambda I+\mu H$, where  $\lambda,\mu$ are two complex numbers. More over, if $T$ is self-adjoint with norm $1$, then $T$ must be $\pm I$or  $\pm H$.
\eT
\bpp
It comes from  Proposition \ref{p1} and \ref{pp} that $H$ commutes with the natural representation of $r\mathrm{Spin}(3)+\RR^3$ over $L^2(\RR^3,\HH)$. Again notice that the representations over $ {}^q\!{H}^\pm_1$ are equivalent to those over $ {}^q\!{H}^\pm_2$, respectively.
By Theorem \ref{t1} (Schur's Lemma),  there exist two complex numbers $\lambda$ and $\mu$ such that
\begin{equation}\label{}
\nonumber  Tf=\lambda f,\forall f\in H^+;\ Tg=\mu g,\forall g\in H^-.
\end{equation}
Moreover we have the follows:
\begin{enumerate}
  \item $\lambda,\mu$ must be real numbers because $T$ is real symmetry;
\item   $|\lambda|=|\mu|=1$ because $||T||=1$.
\end{enumerate}
Then  $\lambda$ , $\mu$
must belong to $\{1,-1\}$. If $T\neq\pm I$. Then $T$ is just the Hilbert transformation when $\lambda=1,\mu=-1$, while $T=-H$ when $\lambda=-1,\mu=1$.
\epp
\brs
We note that as a bi-product the theorem concludes that the operators commuting with the natural representation of $r\mathrm{Spin}(3)+\RR^3$ are automatically bounded operators.
\ers
\section{Symmetry of the Hilbert transformation: case of  $\mathrm{Cl}_2$}
As the same as the discussion in section 5, we obtain the following lemmas.
\bL \label{2de} Let $\chi_\pm(\ul{x})=\frac12(1\pm i\frac{\ul{x}}{|\ul{x}|}),\ \ul{x}\in \RR^{0,2}$.
\begin{enumerate}
  \item $f\in \tilde{\tilde{H}}^\pm_1$ if and only if $f\in L^2(\RR^2,S_2)$ and $f=\chi_\pm f$;
 \item $g\in \tilde{\tilde{H}}^\pm_2$ if and only if $g\in L^2(\RR^2,S_2e_1)$ and $g=\chi_\pm g$.
\end{enumerate}

\eL

\bL \label{2l2s2}
Hilbert space  $L^2(\RR^2,S_2)$ is the direct sum of $\tilde{\tilde{H}}^+_1$ and $\tilde{\tilde{H}}^-_1$.
\eL

\bL \label{2l2s23}
The Hilbert space  $L^2(\RR^2,S_2e_1)$ is the direct sum of $\tilde{\tilde{H}}^+_2$ and $\tilde{\tilde{H}}^-_2$.
\eL

\vskip2mm
\noindent \it Proof of Lemma  \ref{2de}---\ref{2l2s23}.  \rm
It is similar to the proof of  Theorem \ref{de}. For instance, to prove Lemma \ref{2l2s2}  and the first part of Lemma \ref{2de}, we only need to in the proof of Theorem \ref{de}
replace $e_3(1-ie_3)=i(1-ie_3)$  by $e_2(e_1+ie_2+1-ie_1e_2)=i(e_1+ie_2+1-ie_1e_2)$, and $e_3(e_1+ie_2)=-i(e_1+ie_2)$  by $e_2(e_1+ie_2-(1-ie_1e_2))=-i(e_1+ie_2-(1-ie_1e_2))$.

Now we obtain
\bT
\begin{equation}\label{2decH}
  L^2(\RR^2,\mathrm{Cl}_2)=L^2(\RR^2,S_2)\oplus L^2(\RR^2,S_2e_1)=\tilde{\tilde{H}}^+_1\oplus \tilde{\tilde{H}}^-_1\oplus \tilde{\tilde{H}}^+_2\oplus \tilde{\tilde{H}}^-_2.
\end{equation}
\eT

Let $^2\!{H}^\pm_j$ be the Hilbert spaces isometrically isomorphic to $\tilde{\tilde{H}}^\pm_j$, $j=1,2$, as images of the Fourier
transformation mappings, i.e., $^2\!{H}^\pm_j=\mathfrak{F}\tilde{\tilde{H}}^\pm_j$.  Denote
$$
^2\!{H}^+:=\tilde{\tilde{H}}^+_1\oplus\tilde{\tilde{H}}^+_2,  ^2\!{H}^-=\tilde{\tilde{H}}^-_1\oplus\tilde{\tilde{H}}^-_2.
$$
Then we can obtain the following theorems.
\bT Let $H$ be the Hilbert transformation extended to $L^2(\RR^2,\mathrm{Cl}_2)$ based on its definition over $L^2(\RR^2)$ and the relation (\ref{hilbert3}). Then it has two nontrivial maximal invariant
spaces $ {}^2\!{H}^+$ and $ {}^2\!{H}^-$ spelt as
\begin{equation}\label{}
 \nonumber   L^2(\RR^2,S_2)= {}^2\!{H}^+\oplus {}^2\!{H}^-,
\end{equation}
and
\begin{equation}\label{}
  Hf=f,\quad \forall f\in  {}^2\!{H}^+;\quad Hg=-g,\quad \forall g\in  {}^2\!{H}^-.
\end{equation}
\eT

\bT
Let $T$ be a  linear  operator over $L^2(\RR^3,\mathrm{Cl}_2)$. If  $T$ is  commutes with the natural representation of $\mathfrak{G}_2$. Then $T$ must be be of the form $\lambda I+\mu H$, where  $\lambda,\mu$ are two complex numbers. More over, if $T$ is self-adjoint with norm $1$, then $T$ must be $\pm I$or  $\pm H$.
\eT

The proofs of Theorem 6.5 and 6.6 are similar to one of the Theorem 5.9 and \ref{qlast}.
\section{operators with symmetry $\mathfrak{G}$: case of  $\mathrm{Cl}_3$}
 In this section  we continue our discussion on symmetry of the  Hilbert  transformations over $L^2(\RR^3,\mathrm{Cl}_3)$.  If we merely substitute $S_2$ for $S'_2$ in the argument of \S 5, the proof of Lemma \ref{de} is invalid because there do not exist $\lambda,\mu\in \CC$ such that $e_3(e_1+ie_1e_2)=\lambda(e_1+ie_1e_2)$ and $ e_3(1-ie_1e_2)=\mu (1-ie_1e_2)$. Then we can not obtain the counterpart results to the quaternionic case. So we substitute $W_2$ for $S_2$ to revise the arguments.
\subsection{Decomposition of $L^2(\RR^3,\mathrm{Cl}_3)$}
As the same as the case for $\HH$, $L^2(\RR^3,\mathrm{Cl}_3)$ should not been considered as the direct sum of two subspaces in the view of the symmetry of the Hilbert  transformations.
\vspace{2mm}

\noindent We still have the following proposition.
\bP The Hilbert transformation $H$ and $\pi_{rsb},(r,s,b)\in\mathfrak{G},$ are the unitary operators on $L^2(\RR^3,\mathrm{Cl}_3)$. And
\begin{equation}\label{pp2}
  \nonumber (H\pi_{rsb})f=(\pi_{rsb}H)f,\quad f\in L^2(\RR^3)
\end{equation}
if and only if
\begin{equation}\label{}
  \nonumber (H\pi_{rsb})f=(\pi_{rsb}H)f,\quad f\in L^2(\RR^3,\mathrm{Cl}_3).
\end{equation}
\eP

Denote by
$$H'_1:=L^2(\RR^3,\gamma(\f{\ul{y}}{|\ul{y}|}) W_2), $$
$$H'_2:=L^2(\RR^3,\gamma(\f{\ul{y}}{|\ul{y}|}) W_2e_1), $$
$$H'_3:=L^2(\RR^3,\gamma(\f{\ul{y}}{|\ul{y}|}) W_2e_3),$$
$$H'_4:=L^2(\RR^3,\gamma(\f{\ul{y}}{|\ul{y}|}) W_2e_1e_3), $$
$$H'^+:=\oplus_{j=1}^4 H'^+_j,\ H'^-=\oplus_{j=1}^4 H'^-_j.$$
Then we obtain two lemmas.

 \bL \label{de2}
\begin{enumerate}
  \item $f\in H'^+_j$ if and only if $f\in H'_j$ and $f=\chi_+f$, $j=1,2,3,4$;
  \item $f\in H'^-_j$ if and only if $f\in H'_j$ and $f=\chi_-f$,  $j=1,2,3,4$.
  \item $f\in H'^+$ if and only if $f\in L^2(\RR^3,\mathrm{Cl}_3)$ and $f=\chi_+f$;
   \item $f\in H'^-$ if and only if $f\in L^2(\RR^3,\mathrm{Cl}_3)$ and $f=\chi_-f$.
\end{enumerate}

\eL
\bL \label{l2s22}
  $L^2(\RR^3,\mathrm{Cl}_3)=H'^+\bigoplus H'^-=\bigoplus_{j=1}^4 H'^+_j\bigoplus\bigoplus_{j=1}^4 H'^-_j=\bigoplus_{j=1}^4H'_j.$
\eL
\vskip2mm
\noindent  The proofs of Lemma  {\ref{de2} and\ref{l2s22}  are similar to the proof of  Theorem \ref{de}. In the four proof, We only need to
$e_3(1-ie_3)=i(1-ie_3)$ is replaced by, respectively,
 \begin{eqnarray*}
 \nonumber e_3(1-e_3e_1e_2)(1-ie_1e_2)&=&i(1-e_3e_1e_2)(1-ie_1e_2) \\
  \nonumber e_3(1-e_3e_1e_2)(1-ie_1e_2)e_1&=&i(1-e_3e_1e_2)(1-ie_1e_2)e_1\\
  \nonumber e_3(1-e_3e_1e_2)(1-ie_1e_2)e_3&=&i(1-e_3e_1e_2)(1-ie_1e_2)e_3\\
  \nonumber  e_3(1-e_3e_1e_2)(1-ie_1e_2)e_1e_3&=&i(1-e_3e_1e_2)(1-ie_1e_2)e_1e_3,
 \end{eqnarray*}
and $e_3(e_1+ie_2)=-i(e_1+ie_2)$ is replaced by, respectively,
 \begin{eqnarray*}
 \nonumber e_3(1-e_3e_1e_2)((e_1+ie_2)&=&-i(1-e_3e_1e_2)(e_1+ie_2)\\
  \nonumber e_3(1-e_3e_1e_2)((e_1+ie_2)e_1&=&-i(1-e_3e_1e_2)(e_1+ie_2)e_1 \\
   \nonumber e_3(1-e_3e_1e_2)((e_1+ie_2)e_3&=&-i(1-e_3e_1e_2)(e_1+ie_2)e_3
    \end{eqnarray*}
and $e_3(1-e_3e_1e_2)((e_1+ie_2)e_1e_3=-i(1-e_3e_1e_2)(e_1+ie_2)e_1e_3$.

\subsection{Decomposition of $\pi$ } The Fourier theory of Clifford algebra-valued functions is the same as that for quaternion-valued functions.
 Fourier  transformation   acts as an isometrical automorphism between $H'_j,j=1,2,3,4$.
Let ${H}^\pm_j$ be the Hilbert spaces isometrically isomorphic to $H'^\pm_j$, $j=1,2,3,4$, by Fourier
transformation, respectively, i.e., ${H}^\pm_j=\mathfrak{F}H'^\pm_j$. Still denote by $H_j=H'_j, j=1,2,3,4$. Then $H_j$ have another direct sum decompositions
\begin{equation}\label{decompij}
   H_j={H}^+_j\oplus {H}^-_j,\ j=1,2,3,4.
\end{equation}
 Moreover, we obtain that the natural   representation  ${\pi}^\pm$ over ${H}^\pm_j$  are equivalent to
$\tilde{\pi}'^\pm$ respectively through the Fourier transformation.

Denote
\begin{equation}\label{find it}
\nonumber  H^+=\oplus^4_{j=1}{H}^+_j,\ H^-=\oplus^4_{j=1}{H}^-_j.
\end{equation}

\bT\label{m12}
Let $\pi$ be the natural unitary
representation  of $\mathfrak{G}_3$ over $L^2(\RR^3,\mathrm{CL}_3)$.
Then  $\pi$ is decomposed into the direct sum of sub-representations  ${\pi}^\pm$ over  ${H}^\pm$.
And ${\pi}^\pm$ over ${H}^\pm$ are decomposed into the direct sum of the irreducible  sub-representations over ${H}^\pm_j,j=1,2,3,4$, respectively.
\eT

By lemma 7.2 and 7.3, Theorem \ref{m12} can be proved as the same as that of Theorem \ref{ml}.

\bT\label{last2} Let $H$ be the Hilbert transformation over $L^2(\RR^3,\mathrm{CL}_3)$ extended from its form over $L^2(\RR^3)$ in accordance with (\ref{hilbert3}). Then it has two nontrivial maximal invariant
subspaces ${H}^+$ and ${H}^-$ such that
\begin{equation}\label{}
 \nonumber   L^2(\RR^n,\mathrm{CL}_3)={H}^+\oplus {H}^-,
\end{equation}
and
\begin{equation}\label{h-charact2}
  Hf=f,\quad \forall f\in H^+;\quad Hg=-g,\quad \forall g\in H^-.
\end{equation}
\eT
\bpp  The proof of Theorem \ref{last2} is the same as that of Theorem \ref{last}.
 \epp
 \vspace{2mm}

\noindent Again by  (\ref{h-charact2}) and the Plemelj formulas (Theorem \ref{Plemelj theorem}), we get the following corollary.

 \bC $H^\pm$  are the upper and lower monogenic function Hardy spaces on $\RR^{1,3}_\pm$ respectively.\eC

\vspace{3mm}

\noindent We end this section by establishing the inverse of Theorem \ref{last2}.
\bT
Let $T$ be a  linear operator mapping $L^2(\RR^3,\mathrm{CL}_3)$ into itself. If  $T$ is  commutes with the natural representation of $\mathfrak{G}$. Then $T$ must be of the form $\lambda I+\mu H$, where  $\lambda,\mu$ are two complex numbers. If $T$ is self-adjoint with norm $1$, then $T$ must be $\pm I$or  $\pm H$.
\eT
\bpp
The proof is similar to Theorem \ref{qlast}.
\epp
\section{ More about the symmetry of the Hilbert transforation}

At last, we assert that $\mathfrak{G}$ gives rise to the minimal symmetry of the Hilbert transformation $H$ in some sense.

Let $W$ be a connected subgroup of the automorphism group of $\RR^3$. By the positive solution of the 11th problem of Hilbert, we may suppose that $W$ is a Lie group. $W$ has a natural representation over the $L^2(\RR^3,\mathrm{CL}_3)$: $\tilde{\rho}_w f(\ul{x})=f(U^{-1}\ul{x})$ for $w\in W$. But $H^\pm$ are not invariant under $\tilde{\rho}_w $ even if $w$ is a M\"obius transformation, \cite{el}. So we suppose that there exists the intertwining unitary representation $h(w,x)$ of $W$ over $\mathrm{CL}_3$ such that
\begin{equation}\label{lie}
  \rho(w)(f(\ul{x}))=h(w,\ul{x})|J(w,\ul{x})|^{-1/2}f(w^{-1}\ul{x}),\ \ul{x}\in \RR^3,
\end{equation}
is a representation of $W$ over $L^2(\RR^3,S_2)$, where $J(w,\ul{x})$ is the Jacobian matrix of $w$ at $\ul{x}$.

Notice that it is not necessary for the $h(w,\ul{x})$, $\ul{x}\in \RR^3$,  to be the faithful  representation.

Now suppose that $H^\pm$ are invariant under $\rho$. And $\rho$ commutes with the Hilbert transformation.  Then we  declare that every element of $W$ can be extended to an automorphism of the upper half space. Let $K_{a_0+\ul{a}}(x_0+\ul{x})=\f{a_0+\ul{a}-(x_0+\ul{x})}{|a_0+\ul{a}-(x_0+\ul{x}|^{4}}$, $a_0<0$. Then $K_{a_0+\ul{a}}(x_0+\ul{x})\in H^+$. For $\varphi\in W$ we have
\begin{equation}\label{com-proof1}
  c_n\int_{\RR^3} \f{\ul{y}-\ul{x}}{|y-x|^{4}}(\rho(\varphi)K_{a_0+\ul{a}})(\ul{y})dy
=\rho(\varphi)(c_n\int_{\RR^n}\f{\ul{y}-\ul{x}}{|y-x|^{4}}K_{a_0+\ul{a}}(\ul{y})dy).
\end{equation}
Define $\tilde{K}_{a_0+\ul{a}}(x_0+\ul{x})$  by
\begin{equation}\label{com-proof2}
  \tilde{K}_{a_0+\ul{a}}(x_0+\ul{x})=c_n\int_{\RR^3} \f{\ul{y}-(x_0+\ul{x})}{|\ul{y}-(x_0+\ul{x})|^{4}}(\rho(\varphi)K_{a_0+\ul{a}})(\ul{y})dy.
\end{equation}
Then $\tilde{K}_{a_0+\ul{a}}(x_0+\ul{x})$ is also a function in $H^+$.

Denote by $M_{a_0+\ul{a}}$ the image of $\RR^3$ under the map of $\tilde{K}_{a_0+\ul{a}}(x_0+\ul{x})$. And denote by $\Omega_{a_0+\ul{a}}$ the domain with its boundary $M_{a_0+\ul{a}}$. Remember that the regular function $\tilde{K}_{a_0+\ul{a}}(x_0+\ul{x})$ is an open mapping. We obtain that $\tilde{K}_{a_0+\ul{a}}(x_0+\ul{x})$ maps the upper half space onto $\Omega_{a_0+\ul{a}}$. On the other hand, $K_{a_0+\ul{a}}(x_0+\ul{x})$  is an analytic isomorphism of the half space to $\Omega_{a_0+\ul{a}}$. Thus for $x_0+\ul{x}$ in the upper half space there exists only one point $\xi^{a_0}_0+\ul{\xi^{a_0}}$ in the upper half space such that
\begin{equation}\label{com-proof3}
\tilde{K}_{a_0+\ul{a}}(x_0+\ul{x})=K_{a_0+\ul{a}}(\xi^{a_0}_0+\ul{\xi^{a_0}}).
\end{equation}
Denote by $\xi_0+\ul{\xi}=\lim_{a_0\rightarrow 0^-}\xi^{a_0}_0+\ul{\xi^{a_0}}$. We have
\begin{equation}\label{com-proof4}
\tilde{K}_{\ul{a}}(x_0+\ul{x})=K_{\ul{a}}(\xi_0+\ul{\xi}).
\end{equation}

Again recall that the Hilbert transformation commutes with the translation. We obtain that $\xi_0+\ul{\xi}$ is independent of $\ul{a}$.
So we can define a map $\tilde{\varphi}$ from the upper half space to itself as
\begin{equation}\label{com-proof5}
  \tilde{\varphi}(x_0+\ul{x})=\xi_0+\ul{\xi}.
\end{equation}

It is trivial to check that $\tilde{\varphi}$ is regular since both $\tilde{K}_{\ul{a}}$ and  $K_{\ul{a}}$ are regular. Moreover, we have the following lemma.

\bL \label{1-1} $\tilde{\varphi}$ is a one-to-one mapping from the upper half space to itself.\eL
\bpp Denote by $deg({K}_{\ul{a}},M)$ and $deg(\tilde{K}_{\ul{a}},M)$ the topology degrees of  $K_{\ul{a}}$ and $\tilde{K}_{\ul{a}}$, respectively.
Then
\begin{equation}\label{}
  \nonumber deg(\tilde{K}_{\ul{a}},M) =\sum_{\{x_0+\ul{x}:\ \tilde{K}_{ \ul{a}}(x_0+\ul{x})=\xi_0+\ul{\xi}\}}\sgn{\det \tilde{K}_{\ul{a}}(x_0+\ul{x})}
\end{equation}
for any $\xi_0+\ul{\xi}\in\Omega$, where $\det \tilde{K}_{\ul{a}}(x_0+\ul{x})$ is the determinant of $\tilde{K}_{\ul{a}}$ at $x_0+\ul{x}$.
Notice that $K_{\ul{a}}(\xi_0+\ul{\xi})$ is  homotopic to  $\tilde{K}_{\ul{a}}(x_0+\ul{x})$ since $W$ is connected. We obtain that $deg(\tilde{K}_{\ul{a}},M)=deg({K}_{\ul{a}},M)=1$.

Remember that both $\det K_{a_0+\ul{a}}$ and the determinant of the right multiply by constant Clifford number are positive. We obtain that ${\det \tilde{K}_{\ul{a}}(x_0+\ul{x})}$ is always
positive. So there exists one and
only one point $x_0+\ul{x}$ in the upper half space such that
$\tilde{K}_{ \ul{a}}(x_0+\ul{x})=\xi_0+\ul{\xi}$ for any
$\xi_0+\ul{\xi}\in\Omega$. By (\ref{com-proof4}) and
(\ref{com-proof5}) we obtain that $\tilde{\varphi}$ is a one-to-one
mapping. \epp

Now we prove that the natural action of $\tilde{\varphi}$ is invariant over $H^+$. Let $f\in H^+$. Notice that $K_{\ul{y}}{(x_0+\ul{x})}$ is just the Cauchy kernel.  We obtain that
\begin{eqnarray}\label{com-proof6}
\nonumber    \rho(\tilde{\varphi})f(x_0+\ul{x}) &=& \rho(\tilde{\varphi})(c_n\int_{\RR^3}K_{\ul{y}}(x_0+\ul{x})f(\ul{y})dy)\\
   &=& c_n\int_{\RR^3} \rho(\tilde{\varphi})(K_{\ul{y}}(x_0+\ul{x}))f(\ul{y})dy)
\end{eqnarray}
By (\ref{com-proof2}) and (\ref{com-proof4}) we have that $\tilde{K}_{\ul{y}}(x_0+\ul{x})=\rho(\tilde{\varphi})(K_{\ul{y}}(x_0+\ul{x}))$ is analytic on the upper half space. Then (\ref{com-proof6}) implies that $\rho(\tilde{\varphi})f(x_0+\ul{x})$ belongs to $H^+$. In other words the Dirac operator is invariant under $\tilde{\varphi}$.

It is well known that the Dirac operator is conformally invariant. On the other hand a transformation should be conformal if its natural action commutes with the Dirac operator, \cite{cn}.
Moreover $\tilde{\varphi}$ must  belong to
the generalized M\"{o}bius group (consisting of all the M\"{o}bius transformations and rotations), \cite{al}. Then each $\varphi$ is the product of a rotation and a  M\"{o}bius transformation on $\RR^3$ since it is the restriction of $ \tilde{\varphi}$ to $\RR^3$.
Notice that the reverse transformation is not an automorphism  of $\RR^3$. Thus $W$ can be generated by some transformations among rotations, dilations, and translations. It is easy to check that it is reducible for the natural representations of the groups generated by one or two of those transformations. Moreover, by \cite{sw0}, in order to be commutable with the Dirac operator they have to be rotations with a spinor representation. So $W$ is exactly the $r\mathrm{Spin}(n)+\RR^n$  Group.





\noindent\bf Acknowledgements \rm This paper is supported by Macao Government Science and Technology Development Fund, MSAR. Ref. 045/2015/A2; FDCT099£»FDCT079£» NSFC11701597; NSFC11471250.


\begin{thebibliography}{Ah}

\bibitem{abb}
R.~Abreu-Blaya, J.~Bory-Reyes, F.~Brackx, H.~De Schepper and F.~Sommen, {\it Matrix Cauchy and Hilbert transformations in Hermitian quaternionic Clifford analysis}, Complex Variables and Elliptic Equations, 58(2013), 8, 1057---1069.

\bibitem{al}
L.V.~Ahlfors, {\it Old and new in M\"obius groups.} Ann. Acad. Sci. Fenn. Ser. A I Math. 9 (1984), 93---105.

\bibitem{bg}
 R.~Ba\~{n}uelos, G.~ Wang, {\it Sharp inequalities for martingales with applications to the Beurling-Ahlfors and Riesz transformations}, Duke Math. J. 80 (1995), 3, 575---600.
\bibitem{be}
S.R.~Bell, {\it The Cauchy transform, potential theory, and conformal mapping}, Studies in Advanced Mathematics, CRC Press, Boca Raton, FL, 1992.
\bibitem{bds}
F.~ Brackx,  Richard Delanghe, F.~Sommen, {\it Clifford analysis}. Research Notes in Mathematics, 76. Pitman (Advanced Publishing Program), Boston, MA, 1982.
\bibitem{CZ} A.P. Calder\'on and A. Zygmund, {\it On the existence of certain singular integrals,} Acta Mathematica, {\bf 88}, 1, (1), pp85-139, 1952.
\bibitem{car}
 A.~Carbery, S.~Wainger,J.~Wright,  {\it Double Hilbert transformations along polynomial surfaces in $\RR^3$}, Duke Math. J.
   101(2000), 3, 499---513.
\bibitem{ch}
M.~Christ, {\it Hilbert transformations along curves, II: A flat case},
    Duke Math. J., 52(1985),  4, 887---894.
    \bibitem{cn}
    J.~Cnops, {\it An introduction to Dirac operators on manifolds},  Vol. 24, Progress in Mathematical Physics,
 Birkh\"{a}user, Boston, 2002.
 \bibitem{CMcM}
 R. Coifman, A. McIntosh, and Y. Meyer,
    {\it L'int\'{e}grale de Cauchy d\'{e}finit un
op\'{e}rateur born\'{e} sur $L^2$ pour les courbes lipschitziennes,} Ann. of Math. 116 (1982),
361---387.
\bibitem{cq}
    M. Cowling, T. Qian, {\it A class of singular integralson the n-complex unit sphere,}  Scientia Sinica (Series A),  42(1999), 12, 1233---1245.
\bibitem{dlq}
P.~Dang, H.~Liu, T.~Qian, {\it Hilbert Transformation and Representation of \bf ax+b \rm Group}, {\it to appear in } Canadian Mathematical Bulletin.
\bibitem{dss}
R.~Delanghe, F.~Sommen, V.~Sou$\check{c}$ek, {\it
Clifford algebra and spinor-valued functions}, Kluwer Academic,
Netherlands, 1992.

\bibitem{dm}
V.K.~Dobrev, P.~Moylan, {\it Induced Representations and Invariant Integral Operators for $SU (2, 2)$}, Fortschritte der Physik/Progress of Physics 42(1994), 4,  339---392.
\bibitem{el}
S.-L.~Eriksson, H.~Leutwiler,
{\it Hypermonogenic functions and M\"{o}bius transformations. }
Adv. Appl. Clifford Algebras, 11(2001), S2, 67---76.
\bibitem{eh}
A.S.~ Elmabrok,  O.~Hutnik, {\it Induced representations of the affine group and intertwining operators: I. Analytical approach}, J. Phys. A: Math. Theor. 45(2012), 244017.
\bibitem{folland}
G. B.~ Folland, {\it A Course in Abstract Harmonic Analysis}, CRC
Press, 1995.
\bibitem{carnett}
 J.B.~Garnett, {\it Bounded analytic functions}. Revised first edition. Graduate Texts in Mathematics, 236. Springer, New York, 2007.
 \bibitem{gm}
 J.E. ~Gilbert, M.A.M.~Murray, {\it Clifford Algebras and Dirac operators in Harmonic Analysis}, Cambridge University Press, Cambridge, 1991.
\bibitem{gw}
H.~Gluck, F.W.~Warner,   {\it Great circle fibrations of the
three-sphere}, Duke Math. J., 50(1983), 1, 107---132.
\bibitem{ka}
E.~ Kaniuth,  K.F.~ Taylor, {\it Induced representations of locally compact groups},
Cambridge University Press, New York, 2013.
\bibitem{king}
F.W.~King, {\it Hilbert Transformations}, Volume 1,2, Encyclopedia of Mathematics and its Applications 125. Cambridge University Press,  Cambridge, 2009.
\bibitem{kn}
A.W.~Knapp, {\it Representation theory of semisimple
groups. An overview based on examples}, reprint of
the 1986 original, Princeton Landmarks in Mathematics,
Princeton University Press, Princeton, NJ, 2001.
\bibitem{kim}
J.~Kim, {\it Hilbert transformations along curves in the Heisenberg group},  Proceedings London Mathematical Society, 80(1998), 3, 611---642.
\bibitem{kq}
K.I.~Kou, T.~Qian, {\it The Paley-Wiener theorem in $\RR^n$ with the Clifford analysis setting}, J. Funct. Anal., 189(2002), 227---241.
\bibitem{qm}
C.~Li, A.~McIntosh, T.~Qian,{\it
 Clifford algebras, Fourier transformations, and singular Convolution operators on Lipschitz surfaces},  Revista Matematica Iberoamericana, 10(1994), 3, 665---695.
 \bibitem{q-real ball}
  T. Qian, Fourier analysis on starlike Lipschitz surfaces, Journal of Functional Analysis, 183, 370-412 (2001). DOI:10.1006/jfan.2001.3750.
  \bibitem{qy}T. Qian, Y. Yang, {\it Hilbert Transforms on the Sphere With the Clifford Algebra Setting,} Journal of Fourier Analysis and Applications,  (2009) 15: 753-774. DOI: 10.1007/s00041-009-9062-4.

\bibitem{QXYYY} T. Qian, Y. S. Xu, D. Y. Yan, L. X. Yan and B. Yu, Fourier Spectrum Characterization of Hardy Spaces and Applications, Proceedings of the American Mathematical Society, Volume 137, Number 3, March 2009, page 971-980. DOI:10.1090/S0002-9939-08-09544-0.

\bibitem{sa}
S.~Sahi, Elias M.~Stein, {\it Analysis in matrix space and Speh's representations},
invent. math., 101(1990), 1, 379---393.
\bibitem{st}
 E. M.~Stein, {\it Some problems in harmonic analysis}, Harmonic analysis in Euclidean spaces (Proc. Sympos. Pure Math., Williams Coll., Williamstown, Mass., 1978) Proc. Sympos. Pure Math., XXXV, Part, Amer. Math. Soc., Providence, R.I., 1979, pp. 3---20.
\bibitem{stein}
 E.M.~Stein ,{\it Singular Integrals and Differentiability Properties of Functions}, Princeton University, New Jersey, 1970.
  \bibitem{sw0}
  E. M. Stein, G. Weiss, {\it Generalization of the Cauchy-Riemann equations and representations of the rotation group.} Amer. J. Math. 90(1968),1, 163---196.
 \bibitem{sw}
 E. M. Stein, G. Weiss, {\it Introduction to Fourier analysis on Euclidean spaces,} Princeton University, New Jersey, 1971.







\end{thebibliography}
\end{document}